\documentclass[10pt,reqno,12pt,fleqn,tbtags]{amsart}
\usepackage{amssymb}
\usepackage{amsmath}
\usepackage{geometry}
\usepackage{accents}

\newtheorem{theorem}{Theorem}[section]
\newtheorem{lemma}[theorem]{Lemma}
\newtheorem{proposition}[theorem]{Proposition}
\newtheorem{corollary}[theorem]{Corollary}
\theoremstyle{definition}
\newtheorem{definition}[theorem]{Definition}

\theoremstyle{remark}
\newtheorem{remark}[theorem]{Remark}
\numberwithin{equation}{section}

\begin{document}
\title[\textbf{weighted variable exponent Lebesgue spaces}]{\textbf{variable
Muckenhoupt weights with applications in approximation}}
\author[\textbf{Ramazan Akg\"{u}n}]{\textbf{Ramazan Akg\"{u}n}$^{\mathbf{%
\star}}$}
\address{\textbf{ADDRESS: BALIKES\.{I}R UNIVERSITY, FACULTY OF ARTS AND
SCIENCES, DEPARTMENT OF MATHEMATICS, \c{C}A\u{G}I\c{S} YERLE\c{S}KES\.{I},
10145, BALIKES\.{I}R, T\"{U}RK\.{I}YE.}}
\email{\textbf{rakgun@balikesir.edu.tr}}
\keywords{\textbf{Lebesgue spaces with variable exponent, Variable
Muckenhoupt weight, K functional, Averaging Operator, translation of Steklov
mean, best approximation, simultaneous approximation.}\\
$^{\mathbf{\star }}$\textbf{This work was supported by Balikesir University
Research Project 2019/061.}}

\begin{abstract}
Variable Muckenhoupt weights are considered in variable exponent Lebesgue
spaces. Applications are given for polynomial approximation in these spaces.
Boundedness of averaging operator is proved to gain a transference result.
Almost all weighted norm inequalities are obtained using this transference
result. Potential type approximate identities are considered. Translations
of Steklov operator are proved to be bounded in these spaces. K-functional
is a good apparatus for measuring smoothness properties of functions given
in these function classes. Main inequalities of approximation are derived.
\end{abstract}

\maketitle

\section{\textbf{INTRODUCTION}}

Classical Muckenhoupt's $A_{p}$ class of weights have many important
properties for investigations related to weighted norm inequalities
corresponding to several operators. For example, among others, $A_{p}$
weights are precisely the weights for which many singular integral and
maximal operators are bounded in the appropriate weighted space. Properties
of Muckenhoupt's weights can be found in books \cite{garfra85,graMod,stei},
and paper \cite{muc}. Although Muckenhoupt's class $A_{p}$ comprehensively
studied by many mathematicians working on weight theory of operators, there
is still need to wider class of weights for which characterization problems
of (weighted) operators related to variable exponent function spaces. To
overcome such type problems there has been appeared some variable
Muckenhoupt's class $A_{p(\cdot )}$. As in the classical $L_{p}$ case of $%
A_{p}$, there are two type of definition of weighted variable exponent
Lebesgue spaces: (i) weights as multiplier and (ii) weights as a measure.
Corresponding to these two approaches definition of variable Muckenhoupt's
class $A_{p(\cdot )}$ changes suitably. Papers \cite{cufb,cudh} consider the
case (i) and paper \cite{d-h} considers approach (ii). Also, interestingly,
there has been occured (see \cite{sn21}) a mixed type of (i) and (ii). All
these three types of approaches have their own positive and negative
circumstances. Variable Muckenhoupt's class $A_{p(\cdot )}$ of above used in
many applications to solve several contemporary problems of weighted
operator theory.

In this work, we consider approach (ii) and corresponding class of $%
A_{p(\cdot )}$. Main properties, inequalities or equivalences for weights in 
$A_{p(\cdot )}$ has been proved by Diening H\"{a}st\"{o} (\cite{d-h}) in the
case $1<essinf_{x\in \left( -\infty ,\infty \right) }p(x)$ and $esssup_{x\in
\left( -\infty ,\infty \right) }p(x)<\infty $. In this study we will
consider the case $essinf_{x\in \left[ -\pi ,\pi \right] }p(x)\geq 1$ and 
\newline
$esssup_{x\in \left[ -\pi ,\pi \right] }p(x)<\infty $. We give several
inequalities of trigonometrical approximation of periodic functions
belonging to variable exponent Lebesgue spaces with variable Muckenhoupt's
weights. Applications are given for polynomial approximation in these
spaces. Boundedness of averaging operator is proved to gain a transference
result. Almost all weighted norm inequalities are obtained using this
transference result. Potential type approximate identities are considered.
Translations of Steklov operator are proved to be bounded in these spaces. 
\textit{K}-functional is a good apparatus for measuring smoothness
properties of functions given in these function classes. Main inequalities
of approximation are derived.

\section{\textbf{PRELIMINARY DEFINITIONS}}

\subsubsection{\textbf{LOCAL LOG-HOLDER CONTINUITY CONDITION}}

For measurable set $B\subseteq $\textsc{T}$:=\left[ -\pi ,\pi \right] $, let 
$\mathcal{P}\left( B\right) $ be the class of Lebesgue measurable functions $%
p=p\left( x\right) :B\rightarrow \lbrack 1,\infty )$ such that $1\leq
p_{B}^{-}\leq p_{B}^{+}<\infty $ where $p^{-}\left( B\right) :=essinf_{x\in B}p\left( x\right) $ and $p^{+}\left( B\right) :=esssup_{x\in B}p\left( x\right) $. We set $\mathcal{P}:=\mathcal{P}\left( \text{\textsc{T}}\right) $, $p^{-}:=p^{-}\left( \text{%
\textsc{T}}\right) $ and $p^{+}:=p^{+}\left( \text{\textsc{T}}\right) $.

Variable exponent $p\left( \cdot \right) $ is said to be satisfy \textit{%
Local log-H\"{o}lder continuity condition} (\cite{DHHR11}) on $B\subseteq $%
\textsc{T }if there exists a constant $c_{\log }\left( p\right) >0,$
dependent only on $p$, such that%
\begin{equation}
\left\vert p\left( x\right) -p\left( y\right) \right\vert \leq \frac{c_{\log
}\left( p\right) }{\ln \left( e+1/\left\vert x-y\right\vert \right) }\text{
for all }x,y\in B\text{.}  \label{LLHCc}
\end{equation}%
We define $\mathcal{P}^{\log }\left( B\right) :=\left\{ p\in \mathcal{P}%
\left( B\right) :1/p\text{ satisfy (\ref{LLHCc})}\right\} $ and set $%
\mathcal{P}^{\log }:=\mathcal{P}^{\log }\left( \text{\textsc{T}}\right) .$

\subsubsection{\textbf{WEIGHTED VARIABLE EXPONENT LEBESGUE SPACE}}

Suppose that $B\subseteq $\textsc{T} is a measurable set and $\omega
:B\rightarrow \left[ 0,\infty \right] $ is a $2\pi $ periodic weight
function, i.e., $\omega $ is Lebesgue measurable and almost everywhere
(a.e.) positive function on $B$. We define $\omega \left( B\right)
:=\int\nolimits_{B}\omega \left( x\right) dx$. Let $B\left( x,r\right) $ is
a ball in $\mathbb{R}:=\left( -\infty ,\infty \right) $ with center $x$ and
radius $r>0$.

\begin{definition}
We define weighted variable exponent Lebesgue spaces (\cite{d-h}) as a
collection $L_{2\pi ,\omega }^{p\left( \cdot \right) }\left( B\right) $ of $%
2\pi $ periodic measurable functions $f:B\rightarrow \mathbb{R}$ satisfying%
\begin{equation}
\left\Vert f\right\Vert _{B,p\left( \cdot \right) ,\omega }:=\inf \left\{
\alpha >0:\rho _{B,p\left( \cdot \right) ,\omega }\left( \frac{f}{\alpha }%
\right) \leq 1\right\} <\infty  \label{norm}
\end{equation}%
where $\omega $ is a weight function on $B$\textsc{, }$p\in \mathcal{P}$,
measurable $B\subseteq $\textsc{T }and%
\begin{equation*}
\rho _{B,p\left( \cdot \right) ,\omega }\left( f\right)
:=\int\nolimits_{B}\left\vert f\left( x\right) \right\vert ^{p\left(
x\right) }\omega \left( x\right) dx.
\end{equation*}
\end{definition}
We set $L_{2\pi ,\omega }^{p\left( \cdot \right) }$:=$L_{2\pi ,\omega
}^{p\left( \cdot \right) }\left( \text{\textsc{T}}\right) $, $\rho _{p\left(
\cdot \right) ,\omega }\left( f\right) $:=$\rho _{\text{\textsc{T}},p\left(
\cdot \right) ,\omega }\left( f\right) $, $\left\Vert f\right\Vert _{p\left(
\cdot \right) ,\omega }$:=$\left\Vert f\right\Vert _{\text{\textsc{T}}%
,p\left( \cdot \right) ,\omega }$. $L_{2\pi ,\omega }^{p\left( \cdot \right)
}$ is a Banach space (\cite[Theorem 3.2.7]{DHHR11},\cite{Sh12}) with norm (%
\ref{norm}) when $\omega $ is a weight function on \textsc{T} and $p\in 
\mathcal{P}$. We denote $L_{2\pi ,\omega }^{p\left( \cdot \right) }\left(
B\right) =L^{p}\left( B\right) $ when $\omega \equiv 1$, $p\left( x\right)
=p $ is a constant and measurable $B\subseteq $\textsc{T.}

\begin{definition}
Let $p\in \mathcal{P}$, and $\omega $ be a weight on \textsc{T}. We define $%
p^{\prime }\left( x\right) $:=$p\left( x\right) $/($p\left( x\right) $-$1$)
for $p\left( x\right) >1$ and $p^{\prime }\left( x\right) :=\infty $ for $%
p\left( x\right) =1$.
\end{definition}

\subsection{\textbf{VARIABLE MUCKENHOUPT WEIGHT}}

For given $p\in \mathcal{P}$\textsc{,} the class of weights $\omega $
satisfying the variable exponent Muckenhoupt condition (see Diening H\"{a}st%
\"{o} \cite{d-h})%
\begin{equation}
\left[ \omega \right] _{A_{p\left( \cdot \right) }}:=\underset{B\subseteq 
\text{\textsc{T}}}{\sup }\frac{\left\Vert \omega \right\Vert _{L^{1}\left(
B\right) }}{\left\vert B\right\vert ^{p_{B}}}\left\Vert \frac{1}{\omega }%
\right\Vert _{B,\left( p^{\prime }\left( \cdot \right) /p\left( \cdot
\right) \right) }<\infty  \label{Apx}
\end{equation}%
for any $B\subseteq $\textsc{T}, will be denoted by $A_{p\left( \cdot
\right) }:=A_{p\left( \cdot \right) }\left( \text{\textsc{T}}\right) $ where%
\begin{equation*}
\left( p_{B}\right) ^{-1}:=\frac{1}{\left\vert B\right\vert }\int\limits_{B}%
\frac{dx}{p\left( x\right) }.
\end{equation*}%
When $p\left( x\right) =p=$constant, class $A_{p\left( \cdot \right) }$ in (%
\ref{Apx}) turns into the classical Muckenhoupt class $A_{p}$ (\cite{muc}).

Note that different formulations and definitions of variable exponent
Muckenhoupt weights are also known. See for example the papers Cruz-Uribe,
Diening and H\"{a}st\"{o} (\cite{cudh}), Cruz-Uribe, Fiorenza, Neugebauer (%
\cite{cufb}) and Nogayama, Sawano (\cite{sn21} local variable exponent
Muckenhoupt class).

Formulation (\ref{Apx}) is different from definitions given in \cite%
{cufb,cudh,sn21} because the class $A_{p\left( \cdot \right) }$ is
increasing in $p\left( \cdot \right) .$ Namely,%
\begin{equation}
A_{1}\subset A_{p^{-}}\subset A_{p\left( \cdot \right) }\subset
A_{p^{+}}\subset A_{\infty }  \label{inc}
\end{equation}%
where $1<p^{-}\leq p\left( \cdot \right) \leq p^{+}<\infty $, $A_{\infty
}:=\cup _{p\geq 1}A_{p}$ and%
\begin{eqnarray*}
\omega &\in &A_{1}\Leftrightarrow \left[ \mathcal{\omega }\right] _{A_{1}}:=%
\underset{B\subseteq \text{\textsc{T}}}{\sup }\frac{\omega \left( B\right) }{%
\left\vert B\right\vert }ess\sup_{x\in B}\frac{1}{\omega \left( x\right) }%
<\infty , \\
\omega &\in &A_{p}\Leftrightarrow \left[ \mathcal{\omega }\right] _{A_{p}}:=%
\underset{B\subseteq \text{\textsc{T}}}{\sup }\frac{\omega \left( B\right) }{%
\left\vert B\right\vert ^{p}}\left( \int\nolimits_{B}\left( \omega \left(
x\right) \right) ^{-\frac{1}{p-1}}\right) ^{p-1}<\infty .
\end{eqnarray*}%
The property (\ref{inc}) is no longer true for classes considered in \cite%
{cufb,cudh,sn21}.

We define dual weight as $\omega ^{\prime }\left( \cdot \right) :=\left(
\omega \left( \cdot \right) \right) ^{1-p^{\prime }(\cdot )}.$

\subsection{\textbf{AVERAGING OPERATOR AND TRANSFERENCE RESULT}}

\begin{definition}
(\cite[p.96]{HH}) Let $\mathbb{N}:\mathbb{=}\left\{ 1,2,3,\cdots \right\} $
be natural numbers and $\mathbb{N}_{0}:=\mathbb{N\cup }\left\{ 0\right\} $.

(a) A family $Q$ of measurable sets $E\subset \mathbb{R}$ is called locally $%
N$-finite ($N\in \mathbb{N}$) if 
\begin{equation*}
\sum_{E\in Q}\chi _{E}\left( x\right) \leq N
\end{equation*}%
almost everywhere in $\mathbb{R}$ where $\chi _{U}$ is the characteristic
function of the set $U$.

(b) A family $Q$ of open bounded sets $U\subset \mathbb{R}$ is locally $1$%
-finite if and only if the sets $U\in Q$ are pairwise disjoint.

(c) Let $U\subset $\textsc{T} be \ a measurable set and%
\begin{equation*}
A_{U}f:=\frac{1}{\left\vert U\right\vert }\int\limits_{U\cap \text{\textsc{T}%
}}\left\vert f\left( t\right) \right\vert dt.
\end{equation*}

(d) For a family $Q$ of open sets $U\subset $\textsc{T} we define averaging
operator by 
\begin{equation*}
T_{Q}:L_{loc}^{1}\left( T\right) \rightarrow L^{0}\left( T\right) ,
\end{equation*}%
\begin{equation*}
T_{Q}f\left( x\right) :=\sum_{U\in Q}\chi _{U\cap \text{\textsc{T}}}\left(
x\right) A_{U}f,\quad x\in \text{\textsc{T}},
\end{equation*}%
where $L^{0}\left( \text{\textsc{T}}\right) $ is the set of measurable
functions on \textsc{T}.
\end{definition}

\begin{definition}
(i) Throughout this paper symbol $\mathfrak{A}\precapprox \mathcal{B}$ will
mean that there exists a constant C depending only on unimportant parameters
in question such that inequality $\mathfrak{A}\leq $C$\mathcal{B}$ holds.

(ii) We will use symbol $C$ for generic constants that does not depend on
main parameters and changes with placements. By $\mathbb{S}_{i}$ ($%
i=0,1,\cdots ,13$) we denote specific constants that defined particularly in
context. These constants $\mathbb{S}_{i}$ ($i=0,1,\cdots ,13$) are depend on
main parameters of the problem. On the other hand there are some other
specific constants $c_{inc}$, $c_{\log }\left( p\right) $, $\left[ \omega %
\right] _{A_{p\left( \cdot \right) }}$, $\mathfrak{D}$, $\mathfrak{E}$, $%
\mathfrak{S}$ each of them defined in context particularly. We will give
explicit constants in the proofs but these constants are not best constants.

(iii) For a measurable set $A\subset \mathbb{R}$, symbol $\left\vert
A\right\vert $ will represent the Lebesgue measure of $A$.
\end{definition}

\begin{theorem}
\label{onL1}If $p\in \mathcal{P}^{\log }$ and $\omega \in A_{p\left( \cdot
\right) }$, then there exist positive constants depend only on $p,\omega ,$
such that

(i) $L_{2\pi ,\omega }^{p\left( \cdot \right) }\left( B\right) \subset
L^{1}\left( B\right) $ and 
\begin{equation*}
\left\Vert f\right\Vert _{B,1,1}\precapprox \left\Vert f\right\Vert
_{B,p(\cdot ),\omega }
\end{equation*}%
for any subset $B\subseteq $\textsc{T} with $(1/4)<\left\vert B\right\vert
\leq 2;$

(ii) $L_{2\pi ,\omega }^{p\left( \cdot \right) }\subset L^{1}\left( \text{%
\textsc{T}}\right) $ and 
\begin{equation*}
\left\Vert f\right\Vert _{1,1}\precapprox \left\Vert f\right\Vert _{p(\cdot
),\omega }.
\end{equation*}
\end{theorem}

When $\omega \equiv 1$ Theorem \ref{onL1} is well known. See for example 
\cite[Theorem 2.1]{sh08} of Sharapudinov. When $p\left( \cdot \right) $=$p$%
=constant and $\omega \in A_{p}$ Theorem \ref{onL1} was established in Ky 
\cite{Ky97}, Israfilov \cite[Lemma 2]{dmi02}, Kokilashvili Yildirir \cite[In
Proposition 3.3]{kyil10}, and Berkson and Gillespie \cite[Remark 2.12 p. 934]%
{bg03}. In the case $p^{-}>1$ see also Ayd\i n \cite[Proposition 2.1]{ay12}.

\begin{theorem}
\label{Aver}Suppose that $p\in \mathcal{P}^{\log }$, $\omega \in A_{p\left(
\cdot \right) }$, and $f\in L_{2\pi ,\omega }^{p\left( \cdot \right) }$. If $%
Q$ is $1$-finite family of open bounded subsets of $\mathbb{R}$ having
Lebesgue measure $1$, then, the averaging operator $T_{Q}$ is uniformly
bounded in $L_{2\pi ,\omega }^{p\left( \cdot \right) }$, namely,%
\begin{equation*}
\left\Vert T_{Q}f\right\Vert _{p\left( \cdot \right) ,\omega }\precapprox
\left\Vert f\right\Vert _{p\left( \cdot \right) ,\omega }
\end{equation*}%
holds with a positive constant depending only on $p,\omega $.
\end{theorem}

In case of $\omega \equiv 1$ Theorem \ref{Aver} is obtained by Diening
Harjulehto H\"{a}st\"{o} R\r{u}\v{z}i\v{c}ka \cite{DHHR11} and by Harjulehto
H\"{a}st\"{o} \cite{HH}.

\begin{lemma}
\label{dual}(\cite[p.352]{au20})If $p\in \mathcal{P}^{\log }$, $\omega \in
A_{p\left( \cdot \right) }$, then, dual of $L_{2\pi ,\omega }^{p\left( \cdot
\right) }$ is $L_{2\pi ,\omega ^{\prime }}^{p^{\prime }\left( \cdot \right)
} $.
\end{lemma}

\begin{definition}
\label{Aux}(\cite{z76})Let $p\in \mathcal{P}^{\log }$, $\omega \in
A_{p\left( \cdot \right) }$, and $f\in L_{2\pi ,\omega }^{p\left( \cdot
\right) }$. For an $F\in L_{2\pi ,\omega ^{\prime }}^{p^{\prime }\left(
\cdot \right) }\cap C^{\infty }$, $\left\Vert F\right\Vert _{p^{\prime
}\left( \cdot \right) ,\omega ^{\prime }}\leq 1$, we define auxiliary
function%
\begin{equation}
\mathcal{U}_{f,F}\left( u\right) :=\int\nolimits_{\text{\textsc{T}}}f\left(
x+u\right) \left\vert F\left( x\right) \right\vert dx,\quad u\in \text{%
\textsc{T.}}  \label{efef}
\end{equation}
\end{definition}

Auxiliary function $\mathcal{U}_{f,F}\left( \cdot \right) $ of (\ref{efef})
will frequently be used in the rest of the work and we set $\mathbb{S}%
_{0}:=\left\Vert F\right\Vert _{\infty }$ for further references.

\begin{theorem}
\label{Fu}Let $p\in \mathcal{P}^{\log }$, $\omega \in A_{p\left( \cdot
\right) }$, $F\in L_{2\pi ,\omega ^{\prime }}^{p^{\prime }\left( \cdot
\right) }\cap C^{\infty }$, $\left\Vert F\right\Vert _{p^{\prime }\left(
\cdot \right) ,\omega ^{\prime }}\leq 1$, and $f\in L_{2\pi ,\omega
}^{p\left( \cdot \right) }$. Then, the function $\mathcal{U}_{f,F}\left(
u\right) $, $u\in $\textsc{T}, defined in (\ref{efef}) is a uniformly
continuous function on \textsc{T}.
\end{theorem}

Theorem \ref{Fu} was established in the classical Lebesgue spaces by Ditzian 
\cite{z76}.

Transference Result (briefly TR) is the following theorem.

\begin{theorem}
\label{tra}Let $p\in \mathcal{P}^{\log }$, $\omega \in A_{p\left( \cdot
\right) }$ and $f,g\in L_{2\pi ,\omega }^{p\left( \cdot \right) }$. If%
\begin{equation*}
\left\Vert \mathcal{U}_{f,F}\right\Vert _{C\left( \text{\textsc{T}}\right)
}\precapprox \left\Vert \mathcal{U}_{g,F}\right\Vert _{C\left( \text{\textsc{%
T}}\right) },
\end{equation*}%
with an absolute positive constant, then, weighted norm inequality%
\begin{equation*}
\left\Vert f\right\Vert _{p\left( \cdot \right) ,\omega }\precapprox
\left\Vert g\right\Vert _{p\left( \cdot \right) ,\omega }.
\end{equation*}%
holds with a positive constant depend only on $p,\omega $.
\end{theorem}

Theorem \ref{tra} is new also for $\omega \equiv 1$ and/or for $p\left(
\cdot \right) $=$p$=constant and $\omega \in A_{p}$.

\begin{definition}
(a) Let $f\in L_{2\pi ,\omega }^{p\left( \cdot \right) }$, $\lambda >0$, $%
\tau \in \mathbb{R}$, $x\in $\textsc{T }and%
\begin{equation*}
\text{\textsc{S}}_{\lambda ,\tau \text{\ }}f\left( x\right) :=\lambda
\int\nolimits_{x+\tau -1/\left( 2\lambda \right) }^{x+\tau +1/\left(
2\lambda \right) }f\left( t\right) dt\text{,}
\end{equation*}%
(b) If $0<h<\infty $ we define%
\begin{equation}
T_{h}f\left( x\right) :=\frac{1}{h}\int\nolimits_{x}^{x+h}f\left( t\right)
dt.  \label{tehash0}
\end{equation}
\end{definition}

As a corollary of Transference Result (TR) Theorem \ref{tra} we get the
following result.

\begin{theorem}
\label{stekRR}Suppose that $p\in \mathcal{P}^{\log }$, $\omega \in
A_{p\left( \cdot \right) }$, $f\in L_{2\pi ,\omega }^{p\left( \cdot \right)
} $, $\lambda >0$, and $u\in $\textsc{T}. Then, (a) the family of operators%
\begin{equation*}
\left\{ \text{\textsc{S}}_{\lambda ,u}\right\} _{\lambda >0,u\in \text{%
\textsc{T}}}
\end{equation*}%
is uniformly bounded (in $\lambda $ and $u$) in $L_{2\pi ,\omega }^{p\left(
\cdot \right) }$, namely,%
\begin{equation*}
\left\Vert \text{\textsc{S}}_{\lambda ,u\text{ \ }}f\right\Vert _{p\left(
\cdot \right) ,\omega }\precapprox \left\Vert f\right\Vert _{p\left( \cdot
\right) ,\omega }
\end{equation*}%
holds with a positive constant depend only on $p,\omega $.

(b) If $0<h<\infty $, then, we get%
\begin{equation*}
\left\Vert T_{h}f\right\Vert _{p\left( \cdot \right) ,\omega }\precapprox
\left\Vert f\right\Vert _{p\left( \cdot \right) ,\omega }.
\end{equation*}
\end{theorem}

Theorem \ref{stekRR}(a) for $\lambda \geq 1$ and $\left\vert u\right\vert
\leq \pi /\lambda ^{\rho }$ \ ($\rho >0$) was obtained by Sharapudinov \cite[%
Lemma 3.1]{sh08} and for $\lambda \geq 1$ and $u\in \mathbb{R}$ was obtained
by \L enski Szal \cite[Lemma 2]{lsz20}.

\subsection{APPROXIMATE IDENTITIES}

\begin{definition}
\label{convol}Let $f$ and $g$ be two real-valued $2\pi -$periodic measurable
functions on \textsc{T}. We define the convolution $f\ast g$ of $f$ and $g$
by setting $(f\ast g)(x)=\int\nolimits_{\text{\textsc{T}}}f(y)g(x-y)dy$ for $%
x\in $\textsc{T} for which the integral exists in \textsc{T}.
\end{definition}

\begin{definition}
\label{ddd}(\cite{cuf})Let $B$ be a measurable set $B\subseteq $\textsc{T}, $%
\phi \in L^{1}\left( B\right) $ and $\int_{B}\phi \left( t\right) dt=1$. For
each $t>0$ we define $\phi _{t}\left( x\right) =\frac{1}{t}\phi \left( \frac{%
x}{t}\right) $. Such a sequence $\left\{ \phi _{t}\right\} $ will be called
approximate identity. A function%
\begin{equation*}
\tilde{\phi}\left( x\right) =\sup\limits_{\left\vert y\right\vert \geq
\left\vert x\right\vert }\left\vert \phi \left( y\right) \right\vert
\end{equation*}%
will be called radial majorant of $\phi .$ If $\tilde{\phi}\in L^{1}\left(
B\right) $, then, sequence $\left\{ \phi _{t}\right\} $ will be called
potential-type approximate identity.
\end{definition}

Using the same proof of of Corollary 4.6.6 of \cite[p.130]{DHHR11} we can
obtain the following theorem.

\begin{theorem}
\label{bt}(Corollary 4.6.6 of \cite[p.130]{DHHR11}) Suppose $p\in \mathcal{P}%
^{\log }$, $\omega \in A_{p\left( \cdot \right) }$, $f\in L_{2\pi ,\omega
}^{p\left( \cdot \right) }$, and $\phi $ is a potential-type approximate
identity with radial majorant $\tilde{\phi}\in L^{1}$. Then, for any $t>0$,%
\begin{equation*}
\left\Vert f\ast \phi _{t}\right\Vert _{p\left( \cdot \right) ,\omega
}\precapprox \left\Vert \tilde{\phi}\right\Vert _{1}\left\Vert f\right\Vert
_{p\left( \cdot \right) ,\omega }
\end{equation*}%
and%
\begin{equation*}
\underset{t\rightarrow 0}{\lim }\left\Vert f\ast \phi _{t}-f\right\Vert
_{p\left( \cdot \right) }=0
\end{equation*}%
hold with a positive constant depend only on $p,\omega .$
\end{theorem}

\begin{definition}
(i) Let $X$ be a Banach space on \textsc{T}.\textsc{\ }For $r\in \mathbb{N}$%
, we denote by $W_{X}^{r}$ collection of functions $f\in X$ such that $f^{%
\text{ }\left( r-1\right) }$ is absolutely continuous (AC) and $f^{\text{ }%
\left( r\right) }\in X$. We define $\left\Vert f\right\Vert
_{W_{X}^{r}}:=\left\Vert f\right\Vert _{X}+\left\Vert f^{\text{ }\left(
r\right) }\right\Vert _{X}$.

(ii) In particular:\ (a) Let $C\left( \text{\textsc{T}}\right) $ be the
collection of functions $f:$\textsc{T}$\rightarrow \mathbb{R}$ such that $f$
is continuous function on \textsc{T}. In the case $X=C\left( \text{\textsc{T}%
}\right) $ we set $W_{X}^{r}:=W_{C\left( \text{\textsc{T}}\right) }^{r}$.
(b) In the case $p\in \mathcal{P}$, $\omega \in A_{p\left( \cdot \right) }$
and $X=L_{2\pi ,\omega }^{p\left( \cdot \right) }$ we set $%
W_{X}^{r}:=W_{p\left( \cdot \right) ,\omega }^{r}$. $W_{X}^{r}$ becomes a
Banach space with norm $\left\Vert f\right\Vert _{W_{X}^{r}}$ for $X=C\left( 
\text{\textsc{T}}\right) $ and for $X=L_{2\pi ,\omega }^{p\left( \cdot
\right) }$ where $p\in \mathcal{P}$, $\omega \in A_{p\left( \cdot \right) }$.

(c) For $r\in \mathbb{N}$, we define $C^{r}$ consisting of every member $%
f\in C\left( \text{\textsc{T}}\right) $ such that the derivative $f^{\left(
k\right) }$ exists and is continuous on \textsc{T} for $k=1,...,r$. We set $%
C^{\infty }:=\left\{ f\in C^{r}\text{ for any }r\in \mathbb{N}\right\} $.
\end{definition}

\begin{theorem}
\label{pr2}Let $p\in \mathcal{P}$, $\omega \in A_{p\left( \cdot \right) }$,
and $F\in C^{\infty }$ be as in (\ref{efef}).

(a) If $r\in \mathbb{N}$ and $f\in W_{p\left( \cdot \right) ,\omega }^{r}$,
then, $\mathcal{U}_{f,F}\in C^{r}$ and $\left( \mathcal{U}_{f,F}\right)
^{\left( r\right) }=\left( \mathcal{U}_{f^{\left( r\right) },F}\right) .$

(b) $\mathcal{U}_{g\ast h,F}=\left( \mathcal{U}_{g,F}\right) \ast h$ for $%
g,h\in L_{2\pi ,1}^{1}.$
\end{theorem}

\subsection{\textbf{K-FUNCTIONAL}}

\begin{definition}
Let $X$ be a Banach space on \textsc{T}. For $f\in X$ and $\delta >0$, we
define \textbf{Peetre's K-functional} as%
\begin{equation*}
K_{r}\left( f,\delta ,X,W_{X}^{r}\right) :=\inf\limits_{g\in
W_{X}^{r}}\left\{ \left\Vert f-g\right\Vert _{X}+\delta ^{r}\left\Vert
g^{\left( r\right) }\right\Vert _{X}\right\} .
\end{equation*}%
We set $K_{r}\left( f,\delta ,p\left( \cdot \right) ,\omega \right)
:=K_{r}\left( f,\delta ,L_{2\pi ,\omega }^{p\left( \cdot \right)
},W_{p\left( \cdot \right) ,\omega }^{r}\right) $ for $r\in \mathbb{N}$, $%
p\in \mathcal{P}$, $\delta >0$, $\omega \in A_{p\left( \cdot \right) }$ and $%
f\in L_{2\pi ,\omega }^{p\left( \cdot \right) }$. We also set $K_{r}\left(
f,\delta ,C\right) :=K_{r}\left( f,\delta ,C\left( \text{\textsc{T}}\right)
,W_{C\left( \text{\textsc{T}}\right) }^{r}\right) $ for $r\in \mathbb{N}$, $%
\delta >0$, $f\in C\left( \text{\textsc{T}}\right) $.
\end{definition}

The most important property of \textit{K}-functional is the following result.

\begin{lemma}
\label{kftend}(i) If $r\in \mathbb{N}$, $p\in \mathcal{P}$, $\delta >0$, $%
\omega \in A_{p\left( \cdot \right) }$ and $f\in L_{2\pi ,\omega }^{p\left(
\cdot \right) }$ then%
\begin{equation*}
K_{r}\left( f,\delta ,p\left( \cdot \right) ,\omega \right) \rightarrow
0\quad \text{as }\delta \rightarrow 0.
\end{equation*}%
(ii) If $r\in \mathbb{N}$, $\delta >0$ and $f\in C\left( \text{\textsc{T}}%
\right) $ then%
\begin{equation*}
K_{r}\left( f,\delta ,C\right) \rightarrow 0\quad \text{as }\delta
\rightarrow 0.
\end{equation*}
\end{lemma}

See for example the article Gorbachev Ivanov \cite{git20} for a proof.
Proposition \ref{pr2} (a) and method of proof obtained in \cite{git20} give
conclusions of Lemma \ref{kftend}.

\subsection{\textbf{JACKSON-STECHKIN POLYNOMIAL OPERATOR}}

\begin{definition}
\label{Bappr}Let $p\in \mathcal{P}^{\log }$, $\omega \in A_{p\left( \cdot
\right) }$, and $\mathcal{T}_{n}$ be the class of complex trigonometric
polynomials of degree not greater than $n.$ The best approximation error of
functions $f\in L_{2\pi ,\omega }^{p\left( \cdot \right) }$ by $\mathcal{T}%
_{n}$ is defined by%
\begin{equation*}
E_{n}\left( f\right) _{p\left( \cdot \right) ,\omega }:=\inf \left\{
\left\Vert f-T\right\Vert _{p\left( \cdot \right) ,\omega }:T\in \mathcal{T}%
_{n}\right\} \text{,\quad }n\in \mathbb{N}_{0}.
\end{equation*}
\end{definition}

As an example of approximating polynomial we can give the following
Jackson-Stechkin's polynomial operator.

Let $n,r\in \mathbb{N}$ and (\cite[(2.8) p.204]{devore})%
\begin{equation*}
D_{n,r}f(x):=\frac{1}{\pi }\int\nolimits_{\text{\textsc{T}}}\left[
f(t)+\left( -1\right) ^{r+1}\tilde{\Delta}_{u}^{r}\left( f,t\right) \right] 
\mathcal{J}_{r,n}(u)du\in \mathcal{T}_{n}
\end{equation*}%
be the Jackson-Stechkin's polynomialoperator where $r\leq 2m-2$, $m:=\lfloor 
\frac{n}{r}\rfloor +1$, $\mathcal{J}_{r,n}$ is the Jackson kernel (\cite[%
(2.5) p.203]{devore})%
\begin{equation*}
\mathcal{J}_{r,n}(x):=\frac{1}{\varkappa _{r,n}}\left( \frac{\sin (mx/2)}{%
\sin (x/2)}\right) ^{2r}\text{,}\quad \varkappa _{r,n}:=\frac{1}{\pi }%
\int\nolimits_{\text{\textsc{T}}}\left( \frac{\sin (mt/2)}{\sin (t/2)}%
\right) ^{2r}dt,
\end{equation*}%
for $\tilde{\Delta}_{h}^{r}f:=\left( \tilde{T}_{h}-I\right) ^{r}f$, $\tilde{T%
}_{h}f\left( \cdot \right) =f\left( \cdot +h\right) $ and $\lfloor a\rfloor $%
:=max$\left\{ z:z\leq a\text{ and }z\text{ integer}\right\} $.

It is known that (\cite[p.147]{DzSh})%
\begin{equation*}
\frac{3}{2\sqrt{r}}n^{2r-1}\leq \varkappa _{r,n}\leq \frac{5}{2\sqrt{r}}%
n^{2r-1}\text{, (}n\neq 1\text{).}
\end{equation*}%
Jackson kernels $\mathcal{J}_{k,n}$ satisfy relations (\cite[pp:146-149]%
{DzSh})%
\begin{equation}
\left. 
\begin{array}{c}
\frac{1}{\pi }\int_{\text{\textsc{T}}}\mathcal{J}_{r,n}(u)du=1\text{,}%
\medskip \\ 
\frac{1}{\pi }\int_{\text{\textsc{T}}}\left\vert u\right\vert ^{i}\mathcal{J}%
_{r,n}(u)du\mathbb{<}\frac{1}{n^{i}}\text{, \ for }i\leq 2r-2.\medskip%
\end{array}%
\right\}  \label{JK}
\end{equation}

\begin{definition}
Let $0<h<\infty $, $r\in \mathbb{N}$, $p\in \mathcal{P}^{\log }$, $\omega
\in A_{p\left( \cdot \right) }$, $f\in L_{2\pi ,\omega }^{p\left( \cdot
\right) }$, $\delta \geq 0$.

(a) We define difference operator%
\begin{equation*}
\Delta _{h}^{r}f\left( \cdot \right) :=\left( T_{h}-I\right) ^{r}f\left(
\cdot \right)
\end{equation*}%
where $T_{h}f$ is from (\ref{tehash0}), and $I$ is the identity operator.

(b) We define \textit{modulus of smoothness of order }$r$\textit{\ }as%
\begin{equation*}
\begin{tabular}{ll}
$\Omega _{r}\left( f,\delta \right) _{p\left( \cdot \right) ,\omega
}:=\left\Vert \left( I-T_{\delta }\right) ^{r}f\right\Vert _{p\left( \cdot
\right) ,\omega }$;\quad & $\Omega _{0}\left( f,\delta \right) _{p\left(
\cdot \right) ,\omega }:=\left\Vert f\right\Vert _{p\left( \cdot \right)
,\omega }$;$\text{\quad }\delta >0\text{,}$ \\ 
$\Omega _{r}\left( f,0\right) _{p\left( \cdot \right) ,\omega }:=0:=\Omega
_{0}\left( f,0\right) _{p\left( \cdot \right) ,\omega }.$ & 
\end{tabular}%
\end{equation*}
\end{definition}

The following proposition is immediate from Theorem \ref{stekRR}.

\begin{proposition}
\label{ModProp}If $p\in \mathcal{P}^{\log }$, $\omega \in A_{p\left( \cdot
\right) }$, $r\in \mathbb{N}$ and $f\in L_{2\pi ,\omega }^{p\left( \cdot
\right) }$, then,

$(i)$ there exists a positive constant depend only on $p,\omega $ such that%
\begin{equation*}
\left\Vert \left( I-T_{\delta }\right) ^{r}f\right\Vert _{p\left( \cdot
\right) ,\omega }\precapprox \left\Vert f\right\Vert _{p\left( \cdot \right)
,\omega };
\end{equation*}

$\left( ii\right) $\ $\Omega _{r}\left( \cdot ,\delta \right) _{p\left(
\cdot \right) ,\omega }$\ is non-negative and non-decreasing function of $%
\delta \geq 0$;

$\left( iii\right) $\ $\Omega _{r}\left( f_{1}+f_{2},\cdot \right) _{p\left(
\cdot \right) ,\omega }\leq \Omega _{r}\left( f_{1},\cdot \right) _{p\left(
\cdot \right) ,\omega }+\Omega _{r}\left( f_{2},\cdot \right) _{p\left(
\cdot \right) ,\omega }$.
\end{proposition}

\section{\textbf{MAIN RESULTS}}

\begin{theorem}
\label{TSfr}For every $\alpha ,k,r,n-1\in \mathbb{N}_{0}$, $p\in \mathcal{P}%
^{\log }$, $\omega \in A_{p\left( \cdot \right) }$, and $f\in W_{p\left(
\cdot \right) ,\omega }^{\alpha +r}$, inequalities%
\begin{equation}
E_{n}\left( f^{\text{ }(\alpha )}\right) _{p\left( \cdot \right) ,\omega
}\leq \left\Vert f^{\text{ }(\alpha )}-D_{n,r}(f^{\text{ }(\alpha
)})\right\Vert _{p\left( \cdot \right) ,\omega }\precapprox \frac{1}{n^{r}}%
\left\Vert f^{\left( \alpha +r\right) }\right\Vert _{p\left( \cdot \right)
,\omega }  \label{EqnTUR}
\end{equation}%
and%
\begin{equation}
E_{n}\left( f^{\text{ }(\alpha )}\right) _{p\left( \cdot \right) ,\omega
}n^{r}\precapprox E_{n}\left( f^{\text{ }(\alpha +r)}\right) _{p\left( \cdot
\right) ,\omega }  \label{Eqn11}
\end{equation}%
are hold with positive constants that depend only on $k,r,p,\omega .$
\end{theorem}

For $\omega \equiv 1$ (\ref{EqnTUR}) was proved by Sharapudinov \cite{sh13i}
\ and Volosivets \cite{vol17}. (\ref{Eqn11}) also considered in \cite{vol17}
and \cite{sh15}. In the case $p\left( \cdot \right) $=$p$=constant, $\omega
\in A_{p}$ Theorem \ref{TSfr} is known only for $p>1$ see for example papers 
\cite{Akgeja,ra,ra-dmi,akg11gmj,eah86}.

\begin{theorem}
\label{unifDn}If $p\in \mathcal{P}^{\log }$, $\omega \in A_{p\left( \cdot
\right) }$, $f\in L_{2\pi ,\omega }^{p\left( \cdot \right) }$, then,
sequence of Jackson-Stechkin's operator $\{D_{n,k}f\}_{n,k\in \mathbb{N}}$
is uniformly bounded (in $n,k$) in $L_{2\pi ,\omega }^{p\left( \cdot \right)
}$ since%
\begin{equation*}
\left\Vert D_{n,k}f\right\Vert _{p\left( \cdot \right) ,\omega }\precapprox
\left\Vert f\right\Vert _{p\left( \cdot \right) ,\omega }
\end{equation*}%
holds with a positive constant depend only on $p,\omega $.
\end{theorem}

When $\omega \equiv 1$ and $k=2$ Theorem \ref{unifDn} was proved in \cite%
{sh13i}.

Bernstein-Boas-Nikolskii-Civin inequality (\cite{Ni}).

\begin{theorem}
\label{BBCN}\text{Suppose that }$p\in \mathcal{P}^{\log }$, $\omega \in
A_{p\left( \cdot \right) }$, $\alpha \in \mathbb{N}$, $T\in \mathcal{T}_{n}$%
\ and $0<h<2\pi /n$. Then there exists a positive constant depend only on $%
p,\omega $ such that%
\begin{equation}
\left\Vert T^{(\alpha )}\right\Vert _{p\left( \cdot \right) ,\omega
}\precapprox \frac{n^{\alpha }}{2^{\alpha }\sin ^{\alpha }\left( nh/2\right) 
}\left\Vert \Delta _{h}^{\alpha }T\right\Vert _{p\left( \cdot \right)
,\omega }\text{.}  \label{NiSt}
\end{equation}%
In particular, if $h=\pi /n$, then,%
\begin{equation*}
\left\Vert T^{(\alpha )}\right\Vert _{p\left( \cdot \right) ,\omega
}\precapprox n^{\alpha }\left\Vert \Delta _{\pi /n}^{\alpha }T\right\Vert
_{p\left( \cdot \right) ,\omega }\text{.}
\end{equation*}
\end{theorem}

Theorem \ref{BBCN} is new also for $\omega \equiv 1$ and/or $p\left( \cdot
\right) $=$p$=constant, $\omega \in A_{p}.$

\begin{theorem}
\label{modTur}Let $p\in \mathcal{P}^{\log }$, $\omega \in A_{p\left( \cdot
\right) }$, $n,r\in \mathbb{N}$, $\delta >0$, and $T_{n}\in \mathcal{T}_{n}$%
. If $f\in W_{p\left( \cdot \right) ,\omega }^{r}$ and $g\in L_{2\pi ,\omega
}^{p\left( \cdot \right) }$, then, there exist positive constans that depend
only on $r,\omega $\ and $p$\ such that%
\begin{eqnarray}
\left\Vert \Delta _{\delta }^{r}f\right\Vert _{p\left( \cdot \right) ,\omega
} &\precapprox &\delta ^{r}\left\Vert f^{(r)}\right\Vert _{p\left( \cdot
\right) ,\omega }\text{,}  \label{s1} \\
\frac{1}{n^{r}}\left\Vert T_{n}^{(r)}\right\Vert _{p\left( \cdot \right)
,\omega } &\precapprox &K_{r}\left( T_{n},\frac{\pi }{n},p\left( \cdot
\right) ,\omega \right) _{p\left( \cdot \right) ,\omega }\text{ and}
\label{s2} \\
K_{r}\left( g,\delta ,p\left( \cdot \right) ,\omega \right) &\precapprox
&\left\Vert \Delta _{\delta }^{r}g\right\Vert _{p\left( \cdot \right)
,\omega }\precapprox K_{r}\left( g,\delta ,p\left( \cdot \right) ,\omega
\right) .  \label{s3}
\end{eqnarray}
\end{theorem}

When $\omega \equiv 1$ inequality (\ref{s1}) was established in \cite{vol17}
and \cite{sh13i}. Inequalities (\ref{s2}-\ref{s3}) are new also for $\omega
\equiv 1$ and/or $p\left( \cdot \right) $=$p$=constant, $\omega \in A_{p}$.

\begin{theorem}
\label{TBir}Let $\alpha ,k,r,m,n-1\in \mathbb{N}_{0}$, $p\in \mathcal{P}%
^{\log }$, $\omega \in A_{p\left( \cdot \right) }$, and $f\in W_{p\left(
\cdot \right) ,\omega }^{\alpha +r}$. In this case, there exist positive
constants, depend only on $k,\alpha ,p,\omega ,r$, such that%
\begin{equation}
\left\Vert f^{(\alpha )}-D_{n,r}(f^{(\alpha )})\right\Vert _{p\left( \cdot
\right) ,\omega }\precapprox K_{m}\left( f^{\text{ }(\alpha
)},n^{-1},p\left( \cdot \right) ,\omega \right) ,  \label{b3}
\end{equation}%
\begin{equation}
E_{n}\left( f^{(\alpha )}\right) _{p\left( \cdot \right) ,\omega
}n^{r}\precapprox K_{m}\left( f^{\text{ }(\alpha +r)},n^{-1},p\left( \cdot
\right) ,\omega \right)  \label{Eqn12}
\end{equation}%
and a refinement of (\ref{Eqn12})%
\begin{equation}
\left( \prod\nolimits_{s=1}^{n}E_{s}\left( f^{(\alpha )}\right) _{p\left(
\cdot \right) ,\omega }\right) ^{1/n}n^{r}\precapprox K_{m}\left( f^{\text{ }%
(\alpha +r)},n^{-1},p\left( \cdot \right) ,\omega \right)  \label{IJ1}
\end{equation}%
are hold.
\end{theorem}

Theorem \ref{TBir} is also new for $\omega \equiv 1$ and/or $p\left( \cdot
\right) $=$p$=constant, $\omega \in A_{p}$. For functions in $C($\textsc{T}$%
) $ see \cite{TiNa}.

\begin{theorem}
\label{TTers}If\emph{\ }$p\in \mathcal{P}^{\log }$, $\omega \in A_{p\left(
\cdot \right) }$, $r,n\in \mathbb{N}$\emph{\ }and\emph{\ }$f\in L_{2\pi
,\omega }^{p\left( \cdot \right) }$,\ then%
\begin{equation*}
K_{r}\left( f,n^{-1},p\left( \cdot \right) ,\omega \right) n^{r}\precapprox
\sum\nolimits_{\nu =0}^{n}\left( \nu +1\right) ^{r-1}E_{\nu }\left( f\right)
_{p\left( \cdot \right) ,\omega }
\end{equation*}%
holds with a positive constant depend only on $p,\omega $.
\end{theorem}

When $\omega \equiv 1$ Theorem \ref{TTers} was proved in \cite{vol17}. When $%
p^{-}>1$ see also \cite{akg11gmj}. In case $p\left( \cdot \right) $=$p$%
=constant, $\omega \in A_{p}$ Theorem \ref{TTers} was established in \cite%
{Akgeja}. See also papers \cite{ascs,aosss,ahak,ch12,day,isrgu10,isrYir16}.

\begin{theorem}
\label{Ters}Let $r,k\in \mathbb{N}$, $p\in \mathcal{P}^{\log }$, $f\in
L_{2\pi ,\omega }^{p\left( \cdot \right) }$ and $t\in \left( 0,1/2\right) $.
Then, Marchaud type inequality%
\begin{equation*}
K_{r}\left( f,t,p\left( \cdot \right) ,\omega \right) \precapprox
t^{r}\int\nolimits_{t}^{1}\frac{K_{r+k}\left( f,u,p\left( \cdot \right)
,\omega \right) }{u^{r+1}}du
\end{equation*}%
holds with a positive constant depend only on $r,k,p,\omega .$
\end{theorem}

Theorem \ref{Ters} is also new for $\omega \equiv 1$ and/or $p\left( \cdot
\right) $=$p$=constant, $\omega \in A_{p}$. When $p^{-}>1$ one can see
papers \cite{Akgeja,Ja,Ja1,kc16,kol17,yeydmi10}.

Proof of following theorem is standard. See for example \cite{akg11gmj}.

\begin{theorem}
\label{TersTurv}Suppose that $p\in \mathcal{P}^{\log }$, $\omega \in
A_{p\left( \cdot \right) }$\ and\emph{\ }$f\in L_{2\pi ,\omega }^{p\left(
\cdot \right) }$. If series%
\begin{equation}
\sum\nolimits_{\nu =1}^{\infty }\nu ^{r-1}E_{\nu }\left( f\right) _{p\left(
\cdot \right) ,\omega }  \label{Sers}
\end{equation}%
is convergent for some $r\in \mathbb{N}$, then, $f\in W_{p\left( \cdot
\right) ,\omega }^{r}$\ and%
\begin{equation*}
E_{n}\left( f^{\left( r\right) }\right) _{p\left( \cdot \right) ,\omega
}\precapprox \left( n+1\right) ^{r}E_{n}\left( f\right) _{p\left( \cdot
\right) ,\omega }+\sum\nolimits_{\nu =n+1}^{\infty }\nu ^{r-1}E_{\nu }\left(
f\right) _{p\left( \cdot \right) ,\omega }
\end{equation*}%
holds with a positive constant depend only on $p,\omega .$
\end{theorem}

As a corollary of the last theorem we obtain following result.

\begin{corollary}
\label{corTurvTers}If\emph{\ }$p\in \mathcal{P}^{\log }$, $\omega \in
A_{p\left( \cdot \right) }$\ and\emph{\ }$f\in L_{2\pi ,\omega }^{p\left(
\cdot \right) }$, $k\in \mathbb{N}$\ and series (\ref{Sers}) is convergent
for some positive integer $r$, then, there exists a positive constant depend
only on $r,p,\omega $ such that%
\begin{eqnarray*}
K_{k}\left( f^{\left( r\right) },n^{-1},p\left( \cdot \right) ,\omega
\right) &\precapprox &n^{-k}\sum\nolimits_{\nu =0}^{n}\left( \nu +1\right)
^{k+r-1}E_{\nu }\left( f\right) _{p\left( \cdot \right) ,\omega }+ \\
&&+\sum\limits_{\nu =n+1}^{\infty }\nu ^{r-1}E_{\nu }\left( f\right)
_{p\left( \cdot \right) ,\omega }
\end{eqnarray*}%
holds.
\end{corollary}

Corollary \ref{corTurvTers} was proved in \cite{vol17} and in \cite{ra11u}.
When $p^{-}>1$ see for example \cite{akg11gmj}.

Simultaneous approximation results are given in the following.

\begin{theorem}
\label{simlt}Let $r,s\in \mathbb{N}$, $p\in \mathcal{P}^{\log }$, $\omega
\in A_{p\left( \cdot \right) }$ and $f\in W_{p\left( \cdot \right) ,\omega
}^{\alpha +r}$. Suppose that a $t_{n}^{\ast }\in \mathcal{T}_{n}$ satisfy $%
E_{n}\left( f\right) _{p\left( \cdot \right) ,\omega }=\left\Vert
f-t_{n}^{\ast }\right\Vert _{p\left( \cdot \right) ,\omega }.$ Then there
exist a $T\in \mathcal{T}_{2n}$\emph{\ }($n\in N$)\emph{\ }such that for all 
$k=0,1,\ldots ,r$, there hold%
\begin{equation}
\left\Vert f^{(k)}-\left( t_{n}^{\ast }\right) ^{(k)}\right\Vert _{p\left(
\cdot \right) ,\omega }\precapprox \frac{1}{n^{r-k}}E_{n}\left(
f^{(r)}\right) _{p\left( \cdot \right) ,\omega }\text{\quad and}  \label{fc}
\end{equation}%
\begin{equation}
\left\Vert f^{(k)}-T^{(k)}\right\Vert _{p\left( \cdot \right) ,\omega
}\precapprox \frac{1}{n^{r-k}}\Omega _{s}\left( f^{(r)},\frac{1}{n}\right)
_{p\left( \cdot \right) ,\omega }  \label{uc}
\end{equation}%
with some positive constants depend only on $\alpha ,\beta ,p,\omega $.
\end{theorem}

Theorem \ref{simlt} is also new for $\omega \equiv 1$ and/or $p\left( \cdot
\right) $=$p$=constant, $\omega \in A_{p}$. When $p^{-}>1$ see for example 
\cite{Israfil}.

\begin{definition}
\label{defLip}Let $p\in \mathcal{P}^{\log }$, $\omega \in A_{p\left( \cdot
\right) }$\ and\emph{\ }$f\in L_{2\pi ,\omega }^{p\left( \cdot \right) }$.
For $0<\sigma $\ we set $\bar{r}:=\lfloor \sigma /2\rfloor +1$ and
\end{definition}

$Lip\sigma \left( p\left( \cdot \right) ,\omega \right) $:=$\left\{
f\in L_{2\pi ,\omega }^{p\left( \cdot \right) }:\Omega _{\bar{r}}\left(
f,\delta \right) _{p\left( \cdot \right) ,\omega }\precapprox \delta
^{\sigma }\text{,\quad }\delta >0\right\} $.

\begin{corollary}
\label{corLipBir}Let $p\in \mathcal{P}^{\log }$, $\omega \in A_{p\left(
\cdot \right) }$\ and\emph{\ }$f\in L_{2\pi ,\omega }^{p\left( \cdot \right)
}$. If $0<\sigma $, then the following conditions are equivalent:
\end{corollary}

$%
\begin{tabular}{l}
$\left( a\right) \quad f\in Lip\sigma \left( p\left( \cdot \right)
,\omega \right) $, \\ 
$\left( b\right) $\quad $E_{n}\left( f\right) _{p\left( \cdot \right)
,\omega }\precapprox n^{-\sigma }$,\quad $n\in \mathbb{N}$.%
\end{tabular}%
$

\section{\textbf{AUXILIARY RESULTS}}

Here we will collect some auxiliary definions and results required for
proofs.

\begin{definition}
We denote by $S\left( \text{\textsc{T}}\right) $ the collection of simple
functions on \textsc{T}. We set $S_{0}\left( \text{\textsc{T}}\right)
:=\left\{ f\in S\left( \text{\textsc{T}}\right) :f\text{ \ has a compact
support in \textsc{T}}\right\} $.
\end{definition}

From Corollary 3.2.14 of \cite[p.79]{DHHR11}, Remark 3.11 of \cite[p.14]{d-h}
and proof of Lemma 6.7 of \cite[p.23]{d-h} we have the following proposition.

\begin{proposition}
\label{pr1}(Corollary 3.2.14 of \cite[p.79]{DHHR11}) Let $p\in \mathcal{P}%
^{\log }$, $\omega \in A_{p\left( \cdot \right) }$. Then%
\begin{equation}
\frac{1}{2}\left\Vert f\right\Vert _{p\left( \cdot \right) ,\omega }\leq
\sup_{g\in L_{2\pi ,\omega ^{\prime }}^{p^{\prime }\left( \cdot \right)
}:\left\Vert g\right\Vert _{p^{\prime }\left( \cdot \right) ,\omega ^{\prime
}}\leq 1}\int_{\text{\textsc{T}}}\left\vert f\left( x\right) \right\vert
\left\vert g\left( x\right) \right\vert dx\leq 2\left\Vert f\right\Vert
_{p\left( \cdot \right) ,\omega }  \label{zz}
\end{equation}%
holds for all $f\in L_{2\pi ,\omega }^{p\left( \cdot \right) }$.
Furthermore, the supremum in (\ref{zz}) is unchanged if we replace the
condition $g\in L_{2\pi ,\omega ^{\prime }}^{p^{\prime }\left( \cdot \right)
}$ by $g\in S\left( \text{\textsc{T}}\right) $ or $g\in S_{0}\left( \text{%
\textsc{T}}\right) .$
\end{proposition}

Using Theorem \ref{Aver}, Corollary 4.6.6 of \cite[p.130]{DHHR11} and
Theorem \ref{bt} we have the following proposition.

\begin{proposition}
\label{pr1+}Let $p\in \mathcal{P}^{\log }$, $\omega \in A_{p\left( \cdot
\right) }$. Then%
\begin{equation*}
\frac{1}{12\mathbb{S}_{5}}\left\Vert f\right\Vert _{p\left( \cdot \right)
,\omega }\leq \sup_{g\in L_{2\pi ,\omega ^{\prime }}^{p^{\prime }\left(
\cdot \right) }\cap C^{\infty }:\left\Vert g\right\Vert _{p^{\prime }\left(
\cdot \right) ,\omega ^{\prime }}\leq 1}\int_{\text{\textsc{T}}}\left\vert
f\left( x\right) \right\vert \left\vert g\left( x\right) \right\vert dx\leq
2\left\Vert f\right\Vert _{p\left( \cdot \right) ,\omega }
\end{equation*}%
holds for all $f\in L_{2\pi ,\omega }^{p\left( \cdot \right) }$.
\end{proposition}

Convolution of Definition \ref{convol} exists for every $x\in $\textsc{T}
and measurable function. Furthermore, in the classical Lebesgue spaces,
there holds%
\begin{equation*}
\left\Vert f\ast g\right\Vert _{p,1}\leq \left\Vert f\right\Vert
_{p,1}\left\Vert g\right\Vert _{1,1}\text{.}
\end{equation*}%
If $f$ is continuous (respectively absolutely continuous ($\equiv $AC)) then 
$f\ast g$ is continuous (respectively AC).

Let $\mathbb{Z}$ be the set of integers, $\mathbb{Z}^{\ast }$:=$\left\{ z\in 
\mathbb{Z}\text{:}z\neq 0\right\} $ and $\mathbb{Z}_{n}^{\ast }$:=$\left\{
z\in \mathbb{Z}^{\ast }\text{:}\left\vert z\right\vert \leq n\right\} $.

\begin{definition}
Let $S\left[ f\right] $ be the corresponding complex \textit{Fourier} 
\textit{series} of $f\in L^{1}\left( \text{\textsc{T}}\right) $, i.e.,%
\begin{equation*}
S\left[ f\right] \left( x\right) \backsim \sum\limits_{k=-\infty }^{\infty
}c_{k}\left( f\right) e^{ikx}\text{,\quad }c_{k}\left( f\right) \text{:=}%
\frac{1}{2\pi }\int\limits_{\text{\textsc{T}}}f\left( t\right) e^{ikt}dt%
\text{,\quad }k\in \mathbb{Z}\text{.}
\end{equation*}%
We define $\ S_{n}\left( f\right) :=S_{n}\left( x,f\right)
:=\sum\nolimits_{k=-n}^{n}c_{k}\left( f\right) e^{ikx}$,\quad $n\in \mathbb{N%
}_{0}$.
\end{definition}

\begin{lemma}
(\cite{shEm14})Let $B$ be measurable set $B\subseteq $\textsc{T} and $\omega 
$ be a weight function on $B$. For $p\in \mathcal{P}$, $1\leq p\left(
x\right) \leq q\left( x\right) \leq q_{B}^{+}<\infty $ there holds%
\begin{equation}
\left\Vert f\right\Vert _{B,p\left( \cdot \right) ,\omega }\leq \left(
\omega \left( B\right) +1\right) \left\Vert f\right\Vert _{B,q\left( \cdot
\right) ,\omega }  \label{emb}
\end{equation}%
when the left hand side of (\ref{emb}) is finite.
\end{lemma}

\begin{lemma}
(\cite[p.352, Theorem 3]{au20})For measurable set $B\subseteq $\textsc{T, }$%
p\in \mathcal{P}^{\log }\left( B\right) $ and a weight function $\omega $ on 
$B$, the following H\"{o}lder's inequality%
\begin{equation*}
\left\Vert fg\right\Vert _{B,1,1}\leq 2\left\Vert f\right\Vert _{B,p\left(
\cdot \right) ,\omega }\left\Vert g\right\Vert _{B,p^{\prime }\left( \cdot
\right) ,\omega ^{\prime }}
\end{equation*}%
holds for $f\in L_{2\pi ,\omega }^{p\left( \cdot \right) }\left( B\right) $
and $g\in L_{2\pi ,\omega ^{\prime }}^{p^{\prime }\left( \cdot \right)
}\left( B\right) $ and $1=\frac{1}{p\left( \cdot \right) }+\frac{1}{%
p^{\prime }\left( \cdot \right) }$.
\end{lemma}

From Lemma 3.1 of \cite{d-h} we get the following theorem.

\begin{theorem}
\label{cinc}If $p,q\in \mathcal{P}^{\log }$, $q^{-},p^{-}>1$, $q\leq p$ and $%
\omega \in A_{q\left( \cdot \right) }$, then%
\begin{equation*}
\left[ \omega \right] _{A_{p\left( \cdot \right) }}\leq c_{inc}\left[ \omega %
\right] _{A_{q\left( \cdot \right) }}
\end{equation*}%
holds with%
\begin{equation*}
c_{inc}:=\left\{ 
\begin{tabular}{ll}
$16e^{9\left( c_{\log }p^{-1}+c_{\log }q^{-1}\right) }$ & ; $p,q$
nonconstant, \\ 
$16e^{9c_{\log }p^{-1}}$ & ; only $q$ constant.%
\end{tabular}%
\right.
\end{equation*}
\end{theorem}

Using the same proof of Proposition 4.33 of \cite[p.152]{UF13} we get the
following result.

\begin{theorem}
\label{UF13}Let $1\leq p<\infty $, $\omega $ be a weight on \textsc{T} and $%
f\in L_{2\pi ,\omega }^{p}$. In this case,%
\begin{equation*}
\int\limits_{U\cap \text{\textsc{T}}}\left[ A_{U}(\left\vert f\left(
x\right) \right\vert )\right] ^{p}\omega (x)dx\leq \left[ \omega \right]
_{A_{p}}\int\limits_{U\cap \text{\textsc{T}}}\left\vert f\left( x\right)
\right\vert ^{p}\omega (x)dx\text{ \ \ \ }\forall U\subset \text{\textsc{T}}
\end{equation*}%
if and only if $\omega \in A_{p.}$
\end{theorem}

The following result follows directly from (\ref{Apx}).

\begin{lemma}
\label{weightL1} If $p\in \mathcal{P}^{\log }$ and $\omega \in A_{p\left(
\cdot \right) }$, then $\omega \in L^{1}\left( \text{\textsc{T}}\right) $.
\end{lemma}

\begin{lemma}
\label{cSons}(\cite[Corollary of Theorem 4.1]{vmk-sgs03})Suppose that $p\in 
\mathcal{P}^{\log }$, $\omega \in A_{p\left( \cdot \right) }$, $f\in L_{2\pi
,\omega }^{p\left( \cdot \right) }$. In this case the set $C^{\infty }$ of
infinitely continuously differentiable functions on \textsc{T, }is a dense
subset of $L_{2\pi ,\omega }^{p\left( \cdot \right) }$.
\end{lemma}

Using the same proof of Lemma 3.4 of \cite[p.15]{ahh} we have the following
theorem.

\begin{theorem}
\label{adam2}Let $p\in \mathcal{P}^{\log }$ and $\omega \in A_{p\left( \cdot
\right) }$. If $m\in \mathbb{N}$ satisfies $m>2^{p^{+}}\left[ \omega \right]
_{A_{p^{+}}}$, then%
\begin{equation*}
\int\limits_{\text{\textsc{T}}}\frac{\omega (x)dx}{\left( e+\left\vert
x\right\vert \right) ^{m}}\leq \omega \left( B\left( 0,1\right) \right)
2^{2p^{+}}\left[ \omega \right] _{A_{p^{+}}}^{2}\left( \frac{2^{p^{+}}\left[
\omega \right] _{A_{p^{+}}}}{2^{m}-2^{p^{+}}\left[ \omega \right]
_{A_{p^{+}}}}+1\right) .
\end{equation*}
\end{theorem}

Using Theorem \ref{modTur} and Lemma \ref{kftend} we obtain the following
theorem.

\begin{theorem}
If $p\in \mathcal{P}^{\log }$, $\omega \in A_{p\left( \cdot \right) }$ and $%
f\in L_{2\pi ,\omega }^{p\left( \cdot \right) }$, then,%
\begin{equation*}
\lim\limits_{\delta \rightarrow 0^{+}}\Omega _{r}\left( f,\delta \right)
_{p\left( \cdot \right) ,\omega }=0.
\end{equation*}
\end{theorem}

From 2.11 and Theorem 7.1 of \cite[p:37,135]{douux} we have the following
weak type inequality for Hardy Littlewood's maximal operator $M$.

\begin{theorem}
\label{weak}(\cite[p:37,135]{douux})If $\omega \in A_{1}$, $\lambda >0$ and $%
f\in L_{2\pi ,\omega }^{1}$, then,there exists an absolute constant $%
\mathfrak{D}$ (corresponding to the one dimensional case) such that%
\begin{equation*}
\left\Vert \lambda \chi _{\Upsilon _{\lambda }}\right\Vert _{1,\omega }\leq 
\mathfrak{D}\left\Vert f\right\Vert _{1,\omega }
\end{equation*}%
where $\Upsilon _{\lambda }:=\left\{ x\in \text{\textsc{T}}:Mf\left(
x\right) >\lambda \right\} .$
\end{theorem}

\begin{lemma}
\label{bukun} Let $0<h\leq \delta <\infty $, $p\in \mathcal{P}^{\log }$, $%
\omega \in A_{p\left( \cdot \right) }$ and $f\in L_{2\pi ,\omega }^{p\left(
\cdot \right) }$. Then%
\begin{equation*}
\left\Vert \left( I-T_{h}\right) f\right\Vert _{p\left( \cdot \right)
,\omega }\precapprox \left\Vert \left( I-T_{\delta }\right) f\right\Vert
_{p\left( \cdot \right) ,\omega }
\end{equation*}%
holds.
\end{lemma}

\begin{definition}
\label{redelta}Define, for $f\in L_{2\pi ,1}^{1}$, and $\delta >0,$%
\begin{equation*}
\left( \mathfrak{R}_{\delta }f\right) \left( \cdot \right) :=\frac{2}{\delta 
}\int\nolimits_{\delta /2}^{\delta }\left( \frac{1}{h}\int\nolimits_{0}^{h}f%
\left( \cdot +t\right) dt\right) dh.
\end{equation*}%
Using Minkowski's inequality for integrals we obtain the first item of the
following remark.
\end{definition}

\begin{remark}
\label{rr}(a) For $0<\delta <\infty $, $p\in \mathcal{P}^{\log }$, $\omega
\in A_{p\left( \cdot \right) }$, and\emph{\ }$f\in L_{2\pi ,\omega
}^{p\left( \cdot \right) }$ we have%
\begin{equation*}
\left\Vert \mathfrak{R}_{\delta }f\right\Vert _{p\left( \cdot \right)
,\omega }\precapprox \left\Vert f\right\Vert _{p\left( \cdot \right) ,\omega
}
\end{equation*}%
and, hence, $f-\mathfrak{R}_{\delta }f\in L_{2\pi ,\omega }^{p\left( \cdot
\right) }.$

(b) Set $\mathfrak{R}_{\delta }^{r}f:=\left( \mathfrak{R}_{\delta }f\right)
^{r}$ for $r\in \mathbb{N}$.
\end{remark}

\begin{lemma}
\label{bukunA}Let $0<\delta <\infty $, $p\in \mathcal{P}^{\log }$, $\omega
\in A_{p\left( \cdot \right) }$, and\emph{\ }$f\in L_{2\pi ,\omega
}^{p\left( \cdot \right) }$. Then%
\begin{equation*}
\left\Vert \left( I-\mathfrak{R}_{\delta }\right) f\right\Vert _{p\left(
\cdot \right) ,\omega }\precapprox \left\Vert \left( I-T_{\delta }\right)
f\right\Vert _{p\left( \cdot \right) ,\omega }.
\end{equation*}
\end{lemma}

\begin{remark}
\label{difbil}Note that, the function $\mathfrak{R}_{\delta }f$ is
absolutely continuous and differentiable a.e. on \textsc{T} (see \cite[%
Theorem 6.1]{sh13i}).
\end{remark}

The following lemma is obvious from definitions.

\begin{lemma}
\label{akgunArx}Let $0<\delta <\infty $, $p\in \mathcal{P}^{\log }$, $\omega
\in A_{p\left( \cdot \right) }$, and\emph{\ }$f\in W_{p\left( \cdot \right)
,\omega }^{1}$. Then%
\begin{equation*}
\left( \mathfrak{R}_{\delta }f\right) ^{\prime }=\mathfrak{R}_{\delta
}(f^{\prime })\text{\quad and\quad }\left( T_{\delta }f\right) ^{\prime
}=T_{\delta }(f^{\prime })
\end{equation*}%
a.e. (almost everywhere) on \textsc{T}$.$
\end{lemma}

\begin{lemma}
\label{lm05}Let $0<\delta <\infty $, $p\in \mathcal{P}^{\log }$, $\omega \in
A_{p\left( \cdot \right) }$, and\emph{\ }$f\in L_{2\pi ,\omega }^{p\left(
\cdot \right) }$ be given. Then%
\begin{equation}
\delta \left\Vert (\mathfrak{R}_{\delta }f)^{\prime }\right\Vert _{p\left(
\cdot \right) }\precapprox \left\Vert \left( I-T_{\delta }\right)
f\right\Vert _{p\left( \cdot \right) }  \label{ccc}
\end{equation}%
holds.
\end{lemma}

The following lemma can be proved using induction on $r$.

\begin{lemma}
\label{da}Let $0<\delta <\infty $, $r-1\in \mathbb{N}$, and $f\in
W_{1,1}^{r} $ be given. Then%
\begin{equation*}
\left( \mathfrak{R}_{\delta }^{r}f\right) ^{\left( r\right) }=\left( 
\mathfrak{R}_{\delta }\right) ^{\prime }(\mathfrak{R}_{\delta
}^{r-1}f)^{\left( r-1\right) }.
\end{equation*}
\end{lemma}

\section{\textbf{PROOF OF RESULTS}}

\begin{proof}[\textbf{Proof of Theorem \protect\ref{onL1}}]
(i) Let $B\subset $\textsc{T} be a subset and $(1/4)<\left\vert B\right\vert
\leq 2$. If $p_{B}^{\_}>1$, then $p_{B}>1$ and%
\begin{equation*}
\int\nolimits_{B}\left\vert f(x)\right\vert dx\leq 2\max \left\{ 1,\left(
\rho _{B,p^{\prime }\left( \cdot \right) ,\omega }\left( \omega ^{-1}\right)
\right) ^{1/p^{-}}\right\} \left\Vert f\right\Vert _{B,p(\cdot ),\omega }
\end{equation*}%
\begin{equation*}
=2\max \left\{ 1,\max \left\{ 1,\frac{\left[ \omega \right] _{A_{p\left(
\cdot \right) }}\left\vert B\right\vert ^{p_{B}}}{\omega \left( B\right) }%
\right\} ^{p^{+}/p^{-}}\right\} \left\Vert f\right\Vert _{B,p(\cdot ),\omega
}
\end{equation*}%
\begin{equation*}
\leq 2\max \left\{ 1,\frac{\left[ \omega \right] _{A_{p\left( \cdot \right)
}}}{\omega \left( B\left( 0,1\right) \right) }\left( 1+\pi \right)
^{p+}\right\} ^{p^{+}/p^{-}}\left\Vert f\right\Vert _{B,p(\cdot ),\omega
}\leq \mathbb{S}_{1}\left\Vert f\right\Vert _{B,p(\cdot ),\omega }
\end{equation*}%
holds with%
\begin{equation*}
\mathbb{S}_{1}:=\frac{\left[ \omega \right] _{A_{p\left( \cdot \right)
}}2\left( 1+\pi \right) ^{(p^{+})^{2}}}{\left( \omega \left( B\left(
0,1\right) \right) \right) ^{p^{+}/p^{-}}}.
\end{equation*}

We consider the case $p_{B}^{\_}=1$ and $p^{+}>1.$ In this case we can
decompose $B$ as 
\begin{equation*}
B:=\left( \cup _{j}G_{j}\right) \cup \left( \cup _{\iota }b_{\iota }\right)
\cup \left( \cup _{i}N_{i}\right)
\end{equation*}%
where $G_{j}$, $b_{\iota }$ are subintervals of $B,$ and $N_{i}$ are
singletons $meas(N_{i})=0$ and%
\begin{equation*}
G_{j}:=\left\{ x\in B:p\left( x\right) >1\right\} ,\text{ }b_{\iota
}:=\left\{ x\in B:p(x)=1\right\} ,
\end{equation*}%
\begin{equation*}
\left\vert G_{j}\right\vert \leq 2,\text{\quad }\left\vert b_{\iota
}\right\vert \leq 2.
\end{equation*}%
After then we can proceed as%
\begin{equation*}
\int\nolimits_{B}\left\vert f(x)\right\vert
dx=\sum\nolimits_{j}\int\nolimits_{G_{j}}\left\vert f(x)\right\vert
dx+\sum\nolimits_{\iota }\int\nolimits_{b_{\iota }}\left\vert
f(x)\right\vert dx:=\breve{I}+\ddot{I}.
\end{equation*}%
We estimate $\ddot{I}$ first. As $\left. \frac{p^{\prime }\left( \cdot
\right) }{p\left( \cdot \right) }\right\vert _{b_{\iota }}=\infty $ we get%
\begin{equation*}
\ddot{I}=\sum\limits_{\iota }\int\nolimits_{b_{\iota }}\left\vert
f(x)\right\vert dx=\sum\nolimits_{\iota }\int\nolimits_{b_{\iota
}}\left\vert f(x)\right\vert \omega \left( x\right) \omega \left( x\right)
^{-1}dx
\end{equation*}%
\begin{equation*}
\leq \sum\nolimits_{\iota }\left\Vert f\right\Vert _{b_{\iota },p\left(
\cdot \right) ,\omega }\left\Vert \omega ^{-1}\right\Vert _{b_{\iota
},\infty }=\sum\nolimits_{\iota }\frac{\omega \left( b_{\iota }\right) }{%
\left\vert b_{\iota }\right\vert ^{p_{b_{\iota }}}}\left\Vert \frac{1}{%
\omega }\right\Vert _{b_{\iota },\frac{p^{\prime }\left( \cdot \right) }{%
p\left( \cdot \right) }}\frac{\left\vert b_{\iota }\right\vert ^{p_{b_{\iota
}}}}{\omega \left( b_{\iota }\right) }\left\Vert f\right\Vert _{b_{\iota
},p\left( \cdot \right) ,\omega }
\end{equation*}%
\begin{eqnarray*}
&\leq &\frac{2\left( 1+\pi \right) ^{(p^{+})^{2}}\left[ \omega \right]
_{A_{p\left( \cdot \right) }}}{\left( \omega \left( B\left( 0,1\right)
\right) \right) ^{p^{+}/p^{-}}}\sum\nolimits_{\iota }\left\Vert f\right\Vert
_{b_{\iota },p\left( \cdot \right) ,\omega }\leq \frac{2\left( 1+\pi \right)
^{(p^{+})^{2}}\left[ \omega \right] _{A_{p\left( \cdot \right) }}}{\left(
\omega \left( B\left( 0,1\right) \right) \right) ^{p^{+}/p^{-}}}\left\Vert
f\right\Vert _{B,p\left( \cdot \right) ,\omega } \\
&=&\mathbb{S}_{1}\left\Vert f\right\Vert _{B,p\left( \cdot \right) ,\omega }.
\end{eqnarray*}%
For the second expression $\breve{I}$ we find%
\begin{equation*}
\breve{I}=\sum\nolimits_{j}\int\nolimits_{G_{j}}\left\vert f(x)\right\vert
dx\leq \frac{\left[ \omega \right] _{A_{p\left( \cdot \right) }}2\left(
1+\pi \right) ^{(p^{+})^{2}}}{\left( \omega \left( B\left( 0,1\right)
\right) \right) ^{p^{+}/p^{-}}}\sum\limits_{j}\left\Vert f\right\Vert
_{G_{j},p(\cdot ),\omega }
\end{equation*}%
\begin{equation*}
\leq \mathbb{S}_{1}\left\Vert f\right\Vert _{B,p(\cdot ),\omega }.
\end{equation*}%
If $p_{B}^{\_}=1$ and $p^{+}=1$, then $p\equiv 1$. Now result%
\begin{equation*}
\int\nolimits_{B}\left\vert f(x)\right\vert dx\leq \frac{\left[ \omega %
\right] _{A_{1}}\left\vert B\right\vert }{\omega \left( B\right) }\left\Vert
f\right\Vert _{B,1,\omega }
\end{equation*}%
is known from Remark 2.11 of \cite[p.934]{bg03}.

Since%
\begin{equation*}
\dfrac{\left\vert B\right\vert }{\omega \left( B\right) }\text{=}\dfrac{%
\left\vert B\right\vert ^{p_{B}}}{\omega \left( B\right) }\leq \dfrac{\left(
4\mathfrak{D}\right) ^{p^{+}}\left( 1+\pi \right) ^{p^{+}}}{\min \left\{
\left( \omega \left( B\left( 0,1\right) \right) \right) ^{p^{-}/p_{B\left(
0,1\right) }^{+}},\left( \omega \left( B\left( 0,1\right) \right) \right)
^{p^{+}/p_{B\left( 0,1\right) }^{-}}\right\} }\text{:=}\mathbb{S}_{2}
\end{equation*}

we get%
\begin{equation*}
\int\nolimits_{B}\left\vert f(x)\right\vert dx\leq \left[ \omega \right]
_{A_{1}}\mathbb{S}_{2}\left\Vert f\right\Vert _{B,1,\omega }.
\end{equation*}%
As a result, setting $\mathbb{S}_{3}:=\mathbb{S}_{1}\vee \left( \left[
\omega \right] _{A_{1}}\mathbb{S}_{2}\right) $ we have%
\begin{equation*}
\int\nolimits_{B}\left\vert f(x)\right\vert dx\leq \mathbb{S}_{3}\left\Vert
f\right\Vert _{B,1,\omega }
\end{equation*}%
with $a\vee b:=\max \left\{ a,b\right\} .$

(ii) In this case we can decompose \textsc{T:=}$\cup _{j}G_{j}$ where $%
(1/4)<\left\vert G_{j}\right\vert \leq 2.$ Then,%
\begin{equation*}
\int\nolimits_{\text{\textsc{T}}}\left\vert f(x)\right\vert
dx=\sum\nolimits_{j}\int\nolimits_{G_{j}}\left\vert f(x)\right\vert dx\leq
\sum\nolimits_{j}\mathbb{S}_{3}\left\Vert f\right\Vert _{G_{j},p(\cdot
),\omega }
\end{equation*}%
\begin{equation*}
=\mathbb{S}_{3}\sum\nolimits_{j}\left\Vert f\right\Vert _{G_{j},p(\cdot
),\omega }=\mathbb{S}_{3}\left\Vert f\right\Vert _{p(\cdot ),\omega }.
\end{equation*}
\end{proof}

\begin{proof}[\textbf{Proof of Theorem \protect\ref{Aver}}]
Let us consider $f\in L_{2\pi ,\omega }^{p\left( \cdot \right) }$ with $%
\left\Vert f\right\Vert _{p\left( \cdot \right) ,\omega }\leq 1$ and $%
p^{+}>1 $. Suppose that $Q:=\left\{ U:U\text{ open and }\left\vert
U\right\vert =1\right\} $ be a $1$-finite family. We define constant $%
\mathbb{S}_{4}$ as%
\begin{equation*}
\mathbb{S}_{4}:=\mathfrak{E}\left( 2\left( 1+\mathbb{S}_{3}\right) \left(
1+\pi \right) \right) ^{p^{+}}\left( 1+\omega \left( B\left( 0,1\right)
\right) \right) \left( \left( p^{+}\right) ^{\prime }\left[ \omega \right]
_{A_{p^{+}}}\right) ^{\frac{1}{p^{+}-1}}
\end{equation*}%
where absolute constant $\mathfrak{E}>1$ is come from $p^{+}$-Buckley's
univariate estimate of Hardy Littlewood maximal function. Then using
Corollary 2.2.2 of \cite[p.20]{HH} and Theorem \ref{onL1} we obtain%
\begin{equation*}
\rho _{p\left( \cdot \right) ,\omega }\left( \frac{1}{\mathbb{S}_{4}}%
T_{Q}f\right) =\frac{1}{\mathbb{S}_{4}}\int\limits_{\text{\textsc{T}}%
}\left\vert \sum_{U\in Q}\chi _{U\cap \text{\textsc{T}}}\left( x\right)
A_{U}(f)\right\vert ^{p\left( x\right) }\omega (x)dx
\end{equation*}%
\begin{equation*}
\leq \frac{1}{\mathbb{S}_{4}}\sum_{U\in Q}\chi _{U\cap \text{\textsc{T}}%
}\left( x\right) \int\limits_{U\cap \text{\textsc{T}}}\left\vert
A_{U}(f)\right\vert ^{p\left( x\right) }\omega (x)dx
\end{equation*}%
\begin{equation*}
\leq \frac{1}{\mathbb{S}_{4}}\sum_{U\in Q}\chi _{U}\left( x\right)
\int\limits_{U\cap \text{\textsc{T}}}\left( \mathbb{S}_{3}\left\Vert
f\right\Vert _{p\left( \cdot \right) ,\omega }\right) ^{p\left( x\right)
}\omega (x)dx
\end{equation*}%
\begin{equation*}
\leq \frac{1}{\mathbb{S}_{4}}\sum_{U\in Q}\chi _{U}\left( x\right)
\int\limits_{U\cap \text{\textsc{T}}}\left( \left( 1+\mathbb{S}_{3}\right)
^{p^{+}}\left\Vert f\right\Vert _{p\left( \cdot \right) ,\omega }\right)
\omega (x)dx
\end{equation*}%
\begin{equation*}
\leq \frac{\left( 1+\mathbb{S}_{3}\right) ^{p^{+}}}{\mathbb{S}_{4}}%
\sum_{U\in Q}\chi _{U}\left( x\right) \omega (U\cap \text{\textsc{T}})\leq 
\frac{\left( 1+\mathbb{S}_{3}\right) ^{p^{+}}}{\mathbb{S}_{4}}\sum_{U\in
Q}\chi _{U}\left( x\right) \omega (U)
\end{equation*}%
\begin{equation*}
\leq \frac{\left( 1+\mathbb{S}_{3}\right) ^{p^{+}}2^{p^{+}}\mathfrak{E}}{%
\mathbb{S}_{4}}\left( 1+\pi \right) ^{p^{+}}\left( 1+\omega \left( B\left(
0,1\right) \right) \right) \left( \left( p^{+}\right) ^{\prime }\left[
\omega \right] _{A_{p^{+}}}\right) ^{\frac{1}{p^{+}-1}}=1
\end{equation*}%
and hence%
\begin{equation}
\left\Vert T_{Q}f\right\Vert _{p\left( \cdot \right) ,\omega }\leq \mathbb{S}%
_{4}.  \label{aha}
\end{equation}%
General case%
\begin{equation*}
\left\Vert T_{Q}f\right\Vert _{p\left( \cdot \right) ,\omega }\leq \mathbb{S}%
_{4}\left\Vert f\right\Vert _{p(\cdot ),\omega }\text{ \ \ for }f\in L_{2\pi
,\omega }^{p\left( \cdot \right) }
\end{equation*}
can be obtained easily from (\ref{aha}).

If $p^{+}=1,$ then, the following result%
\begin{equation*}
\left\Vert T_{Q}f\right\Vert _{1,\omega }\leq \left[ \omega \right]
_{A_{1}}\left\Vert f\right\Vert _{1,\omega }
\end{equation*}%
follows from Theorem \ref{UF13}:%
\begin{equation*}
\int\limits_{\text{\textsc{T}}}\left\vert \sum_{U\in Q}\chi _{U\cap \text{%
\textsc{T}}}\left( x\right) A_{U}(f)\right\vert \omega (x)dx\leq \sum_{U\in
Q}\chi _{U\cap \text{\textsc{T}}}\left( x\right) \int\limits_{U\cap \text{%
\textsc{T}}}\left\vert A_{U}(f)\right\vert \omega (x)dx
\end{equation*}%
\begin{equation*}
\leq \left[ \omega \right] _{A_{1}}\sum_{U\in Q}\chi _{U}\left( x\right)
\int\nolimits_{U\cap \text{\textsc{T}}}\left\vert f\left( x\right)
\right\vert \omega (x)dx\leq \left[ \omega \right] _{A_{1}}\left\Vert
f\right\Vert _{1,\omega }.
\end{equation*}%
If we set $\mathbb{S}_{5}:=\mathbb{S}_{4}\vee \left[ \omega \right] _{A_{1}}$
and combine results obtained above, then,%
\begin{equation*}
\left\Vert T_{Q}f\right\Vert _{p\left( \cdot \right) ,\omega }\leq \mathbb{S}%
_{5}\left\Vert f\right\Vert _{p(\cdot ),\omega }.
\end{equation*}
\end{proof}

\begin{proof}[\textbf{Proof of Theorem \protect\ref{Fu}}]
(a) Since $C^{\infty }$ is a dense subset of $L_{2\pi ,\omega }^{p\left(
\cdot \right) }$ (see Theorem \ref{cSons}), we consider functions $f\in
C^{\infty }$. For any $\varepsilon >0,$ there exists $\delta :=\delta \left(
\varepsilon \right) >0$ so that $\left\vert f\left( x+u_{1}\right) -f\left(
x+u_{2}\right) \right\vert <\varepsilon $ for any $u_{1},u_{2}\in $\textsc{T}
with $\left\vert u_{1}-u_{2}\right\vert <\delta $. Then, for $F$ of
Definition \ref{Aux}, there holds inequality%
\begin{equation*}
\left\vert \mathcal{U}_{f,F}\left( u_{1}\right) -\mathcal{U}_{f,F}\left(
u_{2}\right) \right\vert \leq \int\nolimits_{\text{\textsc{T}}}\left\vert
f\left( x+u_{1}\right) -f\left( x+u_{2}\right) \right\vert \left\vert
F\left( x\right) \right\vert dx
\end{equation*}%
\begin{equation*}
\leq \max\limits_{x,u_{1},u_{2}\in \text{\textsc{T}}}\left\vert f\left(
x+u_{1}\right) -f\left( x+u_{2}\right) \right\vert \left\Vert F\right\Vert
_{1}\leq \frac{\varepsilon }{\mathbb{S}_{3}}\mathbb{S}_{3}\left\Vert
F\right\Vert _{p^{\prime }\left( \cdot \right) ,\omega ^{\prime }}\leq
\varepsilon
\end{equation*}%
for any $u_{1},u_{2}\in $\textsc{T} with $\left\vert u_{1}-u_{2}\right\vert
<\delta $. Thus conclusion of Theorem \ref{Fu} follows. For the general case 
$f\in L_{2\pi ,\omega }^{p\left( \cdot \right) }$ there exists an $g\in
C^{\infty }$ so that%
\begin{equation*}
\left\Vert f-g\right\Vert _{p\left( \cdot \right) ,\omega }<\frac{\xi }{4%
\mathbb{S}_{3}\mathbb{S}_{0}}
\end{equation*}%
for any $\xi >0$. Therefore%
\begin{equation*}
\left\vert \mathcal{U}_{f,F}\left( u_{1}\right) -\mathcal{U}_{f,F}\left(
u_{2}\right) \right\vert =\left\vert \mathcal{U}_{f,F}\left( u_{1}\right) -%
\mathcal{U}_{g,F}\left( u_{1}\right) \right\vert +\left\vert \mathcal{U}%
_{g,F}\left( u_{1}\right) -\mathcal{U}_{g,F}\left( u_{2}\right) \right\vert +
\end{equation*}%
\begin{eqnarray*}
+\left\vert \mathcal{U}_{g,F}\left( u_{2}\right) -\mathcal{U}_{f,F}\left(
u_{2}\right) \right\vert &=&\left\vert \mathcal{U}_{f-g,F}\left(
u_{1}\right) \right\vert +\frac{\xi }{2}+\left\vert \mathcal{U}%
_{g-f,F}\left( u_{2}\right) \right\vert \\
&\leq &2\mathbb{S}_{3}\mathbb{S}_{0}\left\Vert f-g\right\Vert _{p\left(
\cdot \right) ,\omega }+\frac{\xi }{2}<\xi .
\end{eqnarray*}%
As a result $\mathcal{U}_{f,F}$ is uniformly continuous on \textsc{T}$.$
\end{proof}

\begin{proof}[\textbf{Proof of Theorem \protect\ref{tra}}]
\textbf{\ }Let $0\leq f,g\in L_{2\pi ,\omega }^{p\left( \cdot \right) }$. If 
$\left\Vert g\right\Vert _{p\left( \cdot \right) ,\omega }=0$, then, the
result is obvious. So we may assume that $\left\Vert g\right\Vert _{p\left(
\cdot \right) ,\omega }>0$. In this case, there exists an absolute constant $%
C$ such that%
\begin{align*}
\left\Vert \mathcal{U}_{f,F}\right\Vert _{C\left( \text{\textsc{T}}\right)
}& \leq C\left\Vert \mathcal{U}_{g,F}\right\Vert _{C\left( \text{\textsc{T}}%
\right) }=C\max_{u\in \text{\textsc{T}}}\left\vert \int\nolimits_{\text{%
\textsc{T}}}g\left( x+u\right) \left\vert F\left( x\right) \right\vert
dx\right\vert \\
& =C\mathbb{S}_{0}\left\Vert g\right\Vert _{1}\leq \mathbb{S}_{3}C\mathbb{S}%
_{0}\left\Vert g\right\Vert _{p\left( \cdot \right) ,\omega }.
\end{align*}%
On the other hand, for any $\varepsilon >0$ and appropriately chosen $%
F_{\varepsilon }\in L_{2\pi ,\omega ^{\prime }}^{p^{\prime }\left( \cdot
\right) }$ with 
\begin{equation*}
\int\nolimits_{\text{\textsc{T}}}f\left( x\right) \left\vert F_{\varepsilon
}\left( x\right) \right\vert dx\geq \frac{1}{12\mathbb{S}_{5}}\left\Vert
f\right\Vert _{p\left( \cdot \right) ,\omega }-\varepsilon \text{,\qquad }%
\left\Vert F_{\varepsilon }\right\Vert _{p^{\prime }\left( \cdot \right)
,\omega ^{\prime }}\leq 1\text{,}
\end{equation*}%
(see Proposition \ref{pr2}), one can find%
\begin{equation*}
\left\Vert \mathcal{U}_{f,F}\right\Vert _{C\left( \text{\textsc{T}}\right)
}\geq \left\vert \mathcal{U}_{f,F}\left( 0\right) \right\vert \geq
\int\nolimits_{\text{\textsc{T}}}f\left( x\right) \left\vert F\left(
x\right) \right\vert dx\geq \frac{1}{12\mathbb{S}_{5}}\left\Vert
f\right\Vert _{p\left( \cdot \right) ,\omega }-\varepsilon
\end{equation*}%
In the last inequality we take as $\varepsilon \rightarrow 0+$ and obtain%
\begin{equation*}
\left\Vert \mathcal{U}_{f,F}\right\Vert _{C\left( \text{\textsc{T}}\right)
}\geq \frac{1}{12\mathbb{S}_{5}}\left\Vert f\right\Vert _{p\left( \cdot
\right) ,\omega }\text{.}
\end{equation*}%
Combining these inequalities we get%
\begin{equation*}
\left\Vert f\right\Vert _{p\left( \cdot \right) ,\omega }\leq 12\mathbb{S}%
_{5}\left\Vert \mathcal{U}_{f,F}\right\Vert _{C\left( \text{\textsc{T}}%
\right) }\leq 12\mathbb{S}_{5}C\left\Vert \mathcal{U}_{g,F}\right\Vert
_{C\left( \text{\textsc{T}}\right) }
\end{equation*}%
\begin{equation*}
\leq 12\mathbb{S}_{3}\mathbb{S}_{5}C\mathbb{S}_{0}\left\Vert g\right\Vert
_{p\left( \cdot \right) ,\omega }.
\end{equation*}%
For general $f,g\in L_{2\pi ,\omega }^{p\left( \cdot \right) }$ we obtain%
\begin{equation*}
\left\Vert f\right\Vert _{p\left( \cdot \right) ,\omega }\leq 24\mathbb{S}%
_{3}\mathbb{S}_{5}C\mathbb{S}_{0}\left\Vert g\right\Vert _{p\left( \cdot
\right) ,\omega }.
\end{equation*}
\end{proof}

\begin{proof}[\textbf{Proof of Theorem \protect\ref{stekRR}}]
Proof of (a): Since%
\begin{equation*}
\mathcal{U}_{\text{\textsc{S}}_{\lambda ,\tau }f,F}=\text{\textsc{S}}%
_{\lambda ,\tau }\mathcal{U}_{f,F}
\end{equation*}

we find%
\begin{equation*}
\left\Vert \mathcal{U}_{\text{\textsc{S}}_{\lambda ,\tau }f,F}\right\Vert
_{C\left( \text{\textsc{T}}\right) }=\left\Vert \text{\textsc{S}}_{\lambda
,\tau }(\mathcal{U}_{f,F})\right\Vert _{C\left( \text{\textsc{T}}\right)
}\leq \left\Vert \mathcal{U}_{f,F}\right\Vert _{C\left( \text{\textsc{T}}%
\right) }.
\end{equation*}%
Now from TR we get%
\begin{equation*}
\left\Vert \text{\textsc{S}}_{\lambda ,\tau }f\right\Vert _{p\left( \cdot
\right) ,\omega }\leq 24\mathbb{S}_{3}\mathbb{S}_{5}\mathbb{S}_{0}\left\Vert
f\right\Vert _{p\left( \cdot \right) ,\omega }.
\end{equation*}

Proof of (b) is the same with (a).
\end{proof}

\begin{proof}[\textbf{Proof of Theorem \protect\ref{pr2}}]
(a) follows from Theorem \ref{Fu} and (\ref{efef}).

(b) For any $u\in $\textsc{T}$,$%
\begin{equation*}
\mathcal{U}_{g\ast h,F}=\int\limits_{\text{\textsc{T}}}\left( g\ast h\right)
\left( x+u\right) F\left( x\right) dx=\int\limits_{\text{\textsc{T}}%
}\int\limits_{\text{\textsc{T}}}g(x+u-t)h(t)dtF\left( x\right) dx
\end{equation*}%
\begin{equation*}
=\int\limits_{\text{\textsc{T}}}\int\limits_{\text{\textsc{T}}%
}g(x+u-t)F\left( x\right) dxh(t)dt=\int\limits_{\text{\textsc{T}}%
}F_{g,F}\left( u-t\right) h(t)dt=\left( \mathcal{U}_{g,F}\right) \ast h.
\end{equation*}
\end{proof}

\begin{proof}[\textbf{Proof of Theorem \protect\ref{TSfr}}]
From%
\begin{equation*}
\tilde{\Delta}_{u}^{k}\left( \mathcal{U}_{f^{\left( \alpha \right)
},F},t\right) =\mathcal{U}_{\tilde{\Delta}_{u}^{k}\left( f^{\left( \alpha
\right) },t\right) ,F}
\end{equation*}%
we can write%
\begin{equation*}
\mathcal{U}_{D_{n,r}\left( f^{\left( \alpha \right) }\right) ,F}=D_{n,r}%
\mathcal{U}_{f^{\left( \alpha \right) },F}
\end{equation*}%
and have%
\begin{equation*}
\mathcal{U}_{f^{\left( \alpha \right) }-D_{n,r}\left( f^{\left( \alpha
\right) }\right) ,F}=\mathcal{U}_{f^{\left( \alpha \right) },F}-\mathcal{U}%
_{D_{n,r}\left( f^{\left( \alpha \right) }\right) ,F}=\mathcal{U}_{f^{\left(
\alpha \right) },F}-D_{n,r}\mathcal{U}_{f^{\left( \alpha \right) },F}\text{.}
\end{equation*}%
Define $\omega _{r}\left( f,\delta \right) _{C\left( \text{\textsc{T}}%
\right) }:=\sup_{\left\vert h\right\vert \leq \delta }\left\Vert \left( I-%
\tilde{T}_{h}\right) ^{r}f\right\Vert _{C\left( \text{\textsc{T}}\right) }$.

As a consequence,%
\begin{equation*}
\left\Vert \mathcal{U}_{f^{\left( \alpha \right) }-D_{n,r}\left( f^{\left(
\alpha \right) }\right) ,F}\right\Vert _{C\left( \text{\textsc{T}}\right)
}=\left\Vert \mathcal{U}_{f^{\left( \alpha \right) },F}-D_{n,r}\mathcal{U}%
_{f^{\left( \alpha \right) },F}\right\Vert _{C\left( \text{\textsc{T}}%
\right) }
\end{equation*}%
\begin{equation*}
=\left\Vert \frac{1}{\pi }\int\nolimits_{\text{\textsc{T}}}\tilde{\Delta}%
_{u}^{r}\left( \mathcal{U}_{f^{\left( \alpha \right) },F},t\right) \mathcal{J%
}_{r,n}(u)du\right\Vert _{C\left( \text{\textsc{T}}\right) }
\end{equation*}%
\begin{equation*}
\leq \frac{1}{\pi }\int\limits_{\text{\textsc{T}}}\left\Vert \left( I-\tilde{%
T}_{u}\right) ^{r}\left( \mathcal{U}_{f^{\left( \alpha \right) },F}\right)
\right\Vert _{C\left( \text{\textsc{T}}\right) }\mathcal{J}_{r,n}(u)du
\end{equation*}%
\begin{equation*}
\leq \frac{1}{\pi }\int\limits_{\text{\textsc{T}}}\omega _{r}\left( \mathcal{%
U}_{f^{\left( \alpha \right) },F},\left\vert u\right\vert \right) _{C\left( 
\text{\textsc{T}}\right) }\mathcal{J}_{r,n}(u)du
\end{equation*}%
\begin{equation*}
\leq \frac{1}{\pi }\omega _{r}\left( \mathcal{U}_{f^{\left( \alpha \right)
},F},\frac{1}{n}\right) _{C\left( \text{\textsc{T}}\right) }\int\nolimits_{%
\text{\textsc{T}}}\left( n\left\vert u\right\vert +1\right) ^{r}\mathcal{J}%
_{r,n}(u)du
\end{equation*}%
\begin{equation*}
\leq n^{r}\omega _{r}\left( \mathcal{U}_{f^{\left( \alpha \right) },F},\frac{%
1}{n}\right) _{C\left( \text{\textsc{T}}\right) }\frac{1}{\pi }%
\int\nolimits_{\text{\textsc{T}}}\left( \left\vert u\right\vert +1/n\right)
^{r}\mathcal{J}_{r,n}(u)du
\end{equation*}%
\begin{equation*}
<2^{r+1}\frac{1}{n^{r}}\left\Vert \left( \mathcal{U}_{f^{\left( \alpha
\right) },F}\right) ^{\left( r\right) }\right\Vert _{C\left( \text{\textsc{T}%
}\right) }=2^{r+1}\frac{1}{n^{r}}\left\Vert \mathcal{U}_{f^{\left( \alpha
+r\right) },F}\right\Vert _{C\left( \text{\textsc{T}}\right) }
\end{equation*}%
From this and TR we get righ hand side of (\ref{EqnTUR}):%
\begin{equation*}
\left\Vert f^{\left( \alpha \right) }-D_{n,r}\left( f^{\left( \alpha \right)
}\right) \right\Vert _{p\left( \cdot \right) ,\omega }\leq 2^{r}48\mathbb{S}%
_{3}\mathbb{S}_{5}\mathbb{S}_{0}\frac{1}{n^{r}}\left\Vert f^{\left( \alpha
+r\right) }\right\Vert _{p\left( \cdot \right) ,\omega }.
\end{equation*}%
On the other hand, left hand side of (\ref{EqnTUR}) is obvious.

Let $Q_{n}\in \mathcal{T}_{n}$ be such that%
\begin{equation*}
\left\Vert f^{(\alpha +r)}-Q_{n}\right\Vert _{p\left( \cdot \right) ,\omega
}=E_{n}\left( f^{(\alpha +r)}\right) _{p\left( \cdot \right) ,\omega }\text{%
,\quad }n\in \mathbb{N}.
\end{equation*}%
We suppose%
\begin{equation*}
\phi :=f-I_{\alpha +r}\left[ Q_{n}\right]
\end{equation*}%
where $I_{s}\left[ f\right] $ is the $s$\textit{-th }$\left( s>0\right) $%
\textit{\ integral} of $f$. Then%
\begin{equation*}
\phi ^{(\alpha +r)}=f^{(\alpha +r)}-Q_{n}
\end{equation*}%
and hence%
\begin{equation*}
\left\Vert \phi ^{(\alpha +r)}\right\Vert _{p\left( \cdot \right) ,\omega
}=\left\Vert f^{(\alpha +r)}-Q_{n}\right\Vert _{p\left( \cdot \right)
,\omega }=E_{n}\left( f^{(\alpha +r)}\right) _{p\left( \cdot \right) ,\omega
}.
\end{equation*}%
Therefore we find%
\begin{eqnarray*}
E_{n}\left( \phi ^{(\alpha )}\right) _{p\left( \cdot \right) ,\omega } &\leq
&\frac{48\mathbb{S}_{3}2^{r}\mathbb{S}_{5}\mathbb{S}_{0}}{n^{r}}\left\Vert
\phi ^{(\alpha +r)}\right\Vert _{p\left( \cdot \right) ,\omega } \\
&\leq &\frac{48\mathbb{S}_{3}2^{r}\mathbb{S}_{5}\mathbb{S}_{0}}{n^{r}}%
E_{n}\left( f^{(\alpha +r)}\right) _{p\left( \cdot \right) ,\omega }.
\end{eqnarray*}%
Since%
\begin{equation*}
E_{n}\left( \phi ^{(\alpha )}\right) _{p\left( \cdot \right) ,\omega
}=E_{n}\left( f^{(\alpha )}\right) _{p\left( \cdot \right) ,\omega }
\end{equation*}%
we conclude that (\ref{Eqn11}) holds.
\end{proof}

\begin{proof}[\textbf{Proof of Theorem \protect\ref{unifDn}}]
Let $n,r\in \mathbb{N}$. We get%
\begin{equation*}
\left\Vert \mathcal{U}_{D_{n,r}f,F}\right\Vert _{C\left( \text{\textsc{T}}%
\right) }=\left\Vert D_{n,r}\mathcal{U}_{f,F}\right\Vert _{C\left( \text{%
\textsc{T}}\right) }
\end{equation*}%
\begin{equation*}
=\left\Vert \frac{1}{\pi }\int\nolimits_{\text{\textsc{T}}}\left[ \mathcal{U}%
_{f,F}(t)+\left( -1\right) ^{r+1}\tilde{\Delta}_{u}^{r}\left( \mathcal{U}%
_{f,F},t\right) \right] \mathcal{J}_{r,n}(u)du\right\Vert _{C\left( \text{%
\textsc{T}}\right) }
\end{equation*}%
\begin{equation*}
\leq \frac{1}{\pi }\int\nolimits_{\text{\textsc{T}}}\sum\limits_{v=1}^{r}%
\left\vert \binom{r}{v}\right\vert \left\Vert \mathcal{U}_{f,F}\left( \cdot
+vu\right) \right\Vert _{C\left( \text{\textsc{T}}\right) }\mathcal{J}%
_{r,n}(u)du
\end{equation*}%
\begin{equation*}
\leq \left\Vert \mathcal{U}_{f,F}\right\Vert _{C\left( \text{\textsc{T}}%
\right) }\sum\limits_{v=1}^{r}\left\vert \binom{r}{v}\right\vert \frac{1}{%
\pi }\int\nolimits_{\text{\textsc{T}}}\mathcal{J}_{r,n}(u)du\leq
2^{r}\left\Vert \mathcal{U}_{f,F}\right\Vert _{C\left( \text{\textsc{T}}%
\right) }.
\end{equation*}%
Now, transference result TR gives%
\begin{equation*}
\left\Vert D_{n,r}f\right\Vert _{p\left( \cdot \right) ,\omega }\leq 24%
\mathbb{S}_{3}\mathbb{S}_{5}\mathbb{S}_{0}2^{r}\left\Vert f\right\Vert
_{p\left( \cdot \right) ,\omega }
\end{equation*}
\end{proof}

\begin{proof}[\textbf{Proof of Theorem \protect\ref{BBCN}}]
Let $n\in \mathbb{N}$ and%
\begin{equation*}
T_{n}\left( x\right) =c_{0}+\sum\limits_{\nu \in \mathbb{Z}_{n}^{\ast
}}c_{\nu }e^{i\nu x}\in \mathcal{T}_{n}\text{ (with respect to }x\text{)}
\end{equation*}%
Then%
\begin{equation*}
\mathcal{U}_{T_{n},F}\left( u\right) =c_{0}\tilde{c}_{0}+\sum\limits_{\nu
\in \mathbb{Z}_{n}^{\ast }}c_{\nu }\tilde{c}_{v}e^{i\nu u}\in \mathcal{T}_{n}%
\text{ (with respect to }u\text{)}
\end{equation*}%
Without loss of generality we can take $c_{0}=0$ in this proof. Then%
\begin{equation*}
\left\Vert \mathcal{U}_{T_{n}^{\left( r\right) },F}\right\Vert _{C\left( 
\text{\textsc{T}}\right) }=\left\Vert \left( \mathcal{U}_{T_{n},F}\right)
^{\left( r\right) }\right\Vert _{C\left( \text{\textsc{T}}\right) }\leq
\left( \frac{n}{2\sin \left( nh/2\right) }\right) ^{r}\left\Vert \tilde{%
\Delta}_{h}^{r}\left( \mathcal{U}_{T_{n},F}\right) \right\Vert _{C\left( 
\text{\textsc{T}}\right) }
\end{equation*}%
\begin{equation*}
\leq \mathfrak{S}\left( \frac{n}{2\sin \left( nh/2\right) }\right)
^{r}\left\Vert \Delta _{h}^{r}\left( \mathcal{U}_{T_{n},F}\right)
\right\Vert _{C\left( \text{\textsc{T}}\right) }
\end{equation*}%
\begin{equation*}
\leq \mathfrak{S}\left( \frac{n}{2\sin \left( nh/2\right) }\right)
^{r}\left\Vert \mathcal{U}_{\Delta _{h}^{r}T_{n},F}\right\Vert _{C\left( 
\text{\textsc{T}}\right) }=\left\Vert \mathcal{U}_{\mathfrak{S}\frac{n^{r}}{%
2^{r}\sin ^{r}\left( nh/2\right) }\Delta _{h}^{r}T_{n},F}\right\Vert
_{C\left( \text{\textsc{T}}\right) }
\end{equation*}%
with%
\begin{equation*}
\mathfrak{S:=}\left\{ 
\begin{tabular}{ll}
$2^{r}\left( r^{r}+(34)^{r}\right) $ & ; $r>1$, \\ 
$36$ & ; $r=1$.%
\end{tabular}%
\right.
\end{equation*}

Using TR we get%
\begin{equation*}
\left\Vert T_{n}^{\left( r\right) }\right\Vert _{p(\cdot ),\omega }\leq 
\mathfrak{S}24\mathbb{S}_{3}\mathbb{S}_{5}\mathbb{S}_{0}\frac{n^{r}}{%
2^{r}\sin ^{r}\left( nh/2\right) }\left\Vert \Delta _{h}^{r}T_{n}\right\Vert
_{p(\cdot ),\omega }.
\end{equation*}
\end{proof}

\begin{proof}[\textbf{Proof of Theorem \protect\ref{modTur}}]
Proof of (\ref{s1}): We note that the following inequality is easy to prove
(see \cite[Lemma 3.2]{AG})%
\begin{equation*}
\left\Vert \left( I-T_{\delta }\right) f\right\Vert _{p(\cdot ),\omega }\leq
12\mathbb{S}_{3}\mathbb{S}_{5}\mathbb{S}_{0}\delta \left\Vert f^{\text{ }%
\prime }\right\Vert _{p(\cdot ),\omega }
\end{equation*}%
for $\delta >0$ and $f\in W_{p\left( \cdot \right) ,\omega }^{1}.$Then%
\begin{equation*}
\left\Vert \left( I-T_{\delta }\right) ^{r}f\right\Vert _{p(\cdot ),\omega
}\leq ...\leq (12)^{r}\mathbb{S}_{3}^{r}\mathbb{S}_{5}^{r}\mathbb{S}%
_{0}^{r}\delta ^{r}\left\Vert f^{(r)}\right\Vert _{p(\cdot ),\omega }\text{,}%
\quad \delta >0\text{, \ }\forall f\in W_{p\left( \cdot \right) ,\omega
}^{r}.
\end{equation*}

Proof of (\ref{s2}) is follows from (\ref{NiSt}) and (\ref{s3}).

Proof of (\ref{s3}): We will follow Theorem 1.4 of \cite{ra19}. For $r\in 
\mathbb{N}$ we consider the operator%
\begin{equation*}
\mathcal{A}_{\delta }^{r}\text{:=}I-\left( I-\mathfrak{R}_{\delta
}^{r}\right) ^{r}\text{=}\sum\limits_{j=0}^{r-1}\left( -1\right) ^{r-j+1}%
\binom{r}{j}\mathfrak{R}_{\delta }^{r\left( r-j\right) }\text{.}
\end{equation*}%
From the identity%
\begin{equation*}
I-\mathfrak{R}_{\delta }^{r}=\left( I-\mathfrak{R}_{\delta }\right)
\sum\nolimits_{j=0}^{r-1}\mathfrak{R}_{\delta }^{j}
\end{equation*}%
we find%
\begin{equation*}
\left\Vert \left( I-\mathfrak{R}_{\delta }^{r}\right) g\right\Vert _{p\left(
\cdot \right) ,\omega }\leq \left( \sum\limits_{j=0}^{r-1}\left( 24\mathbb{S}%
_{3}\mathbb{S}_{5}\mathbb{S}_{0}\right) ^{j}\right) \left\Vert \left( I-%
\mathfrak{R}_{\delta }\right) g\right\Vert _{p\left( \cdot \right) ,\omega }
\end{equation*}%
\begin{equation}
\leq \left( \sum\limits_{j=0}^{r-1}\left( 24\mathbb{S}_{3}\mathbb{S}_{5}%
\mathbb{S}_{0}\right) ^{j}\right) 1728\mathbb{S}_{3}\mathbb{S}_{5}\mathbb{S}%
_{0}\left\Vert \left( I-T_{\delta }\right) g\right\Vert _{p\left( \cdot
\right) ,\omega }  \label{fdd}
\end{equation}%
when $0<\delta <\infty $, $p\in \mathcal{P}^{\log }$, $\omega \in A_{p\left(
\cdot \right) }$ and $g\in L_{2\pi ,\omega }^{p\left( \cdot \right) }.$

Using $\left\Vert f-\mathcal{A}_{\delta }^{r}f\right\Vert _{p\left( \cdot
\right) ,\omega }=\left\Vert \left( I-\mathfrak{R}_{\delta }^{r}\right)
^{r}f\right\Vert _{p\left( \cdot \right) ,\omega }$, one obtains%
\begin{equation*}
\left\Vert f-\mathcal{A}_{\delta }^{r}f\right\Vert _{_{p\left( \cdot \right)
,\omega }}=\left\Vert \left( I-\mathfrak{R}_{\delta }^{r}\right)
^{r}f\right\Vert _{p\left( \cdot \right) ,\omega }\leq \cdots 
\end{equation*}%
\begin{equation*}
\cdots \leq \left[ \left( \sum\nolimits_{j=0}^{r-1}\left( 24\mathbb{S}_{3}%
\mathbb{S}_{5}\mathbb{S}_{0}\right) ^{j}\right) 1728\mathbb{S}_{3}\mathbb{S}%
_{5}\mathbb{S}_{0}\right] ^{r}\left\Vert \left( I-T_{\delta }\right)
^{r}f\right\Vert _{p\left( \cdot \right) ,\omega }
\end{equation*}%
\begin{equation*}
=:\mathbb{S}_{6}\left\Vert \left( I-T_{\delta }\right) ^{r}f\right\Vert
_{p\left( \cdot \right) ,\omega }.
\end{equation*}%
On the other hand, using Lemmas \ref{da} and \ref{ccc}, recursively,%
\begin{equation*}
\delta ^{r}\left\Vert \frac{d^{r}}{dx^{r}}\mathfrak{R}_{\delta
}^{r}f\right\Vert _{p\left( \cdot \right) ,\omega }\leq \left( 24\mathbb{S}%
_{3}\mathbb{S}_{5}\mathbb{S}_{0}\right) ^{r}\left\Vert \left( I-T_{\delta
}\right) ^{r}f\right\Vert _{p\left( \cdot \right) ,\omega }.
\end{equation*}%
Thus%
\begin{equation*}
K_{r}\left( f,\delta ;p\left( \cdot \right) ,\omega \right) \leq \left\Vert
f-\mathcal{A}_{\delta }^{r}f\right\Vert _{p\left( \cdot \right) ,\omega
}+\delta ^{r}\left\Vert \left( \mathcal{A}_{\delta }^{r}f\right) ^{\left(
r\right) }\right\Vert _{p\left( \cdot \right) ,\omega }
\end{equation*}%
\begin{equation*}
\leq \mathbb{S}_{6}\left\Vert \left( I-T_{\delta }\right) ^{r}f\right\Vert
_{p\left( \cdot \right) ,\omega }+\sum\nolimits_{j=0}^{r-1}\left\vert \binom{%
r}{j}\right\vert \delta ^{r}\left\Vert \frac{d^{r}}{dx^{r}}\mathfrak{R}%
_{\delta }^{r\left( r-j\right) }f\left( x\right) \right\Vert _{p\left( \cdot
\right) ,\omega }
\end{equation*}%
\begin{equation*}
\leq \mathbb{S}_{6}\left\Vert \left( I-T_{\delta }\right) ^{r}f\right\Vert
_{p\left( \cdot \right) ,\omega }+
\end{equation*}%
\begin{equation*}
+\left( \sum\nolimits_{j=0}^{r-1}\left\vert \binom{r}{j}\right\vert \left( 24%
\mathbb{S}_{3}\mathbb{S}_{5}\mathbb{S}_{0}\right) ^{r-j}\right) \left( 24%
\mathbb{S}_{3}\mathbb{S}_{5}\mathbb{S}_{0}\right) ^{r}\left\Vert \left(
I-T_{\delta }\right) ^{r}f\right\Vert _{p\left( \cdot \right) ,\omega }
\end{equation*}%
\begin{equation*}
\mathbb{=}\max \left\{ \mathbb{S}_{6},\mathbb{S}_{7}\right\} \left\Vert
\left( I-T_{\delta }\right) ^{r}f\right\Vert _{p\left( \cdot \right) ,\omega
}
\end{equation*}%
with $\mathbb{S}_{7}:=\left( \sum\nolimits_{j=0}^{r-1}\left\vert \binom{r}{j}%
\right\vert \left( 24\mathbb{S}_{3}\mathbb{S}_{5}\mathbb{S}_{0}\right)
^{r-j}\right) \left( 24\mathbb{S}_{3}\mathbb{S}_{5}\mathbb{S}_{0}\right)
^{r}.$

For the reverse of the last inequality, when $g\in W_{p\left( \cdot \right)
}^{r}$,%
\begin{equation*}
\Omega _{r}\left( f,\delta \right) _{p\left( \cdot \right) ,\omega }\leq
\left( 1+12\mathbb{S}_{3}\mathbb{S}_{5}\mathbb{S}_{0}\right) ^{r}\left\Vert
f-g\right\Vert _{p\left( \cdot \right) ,\omega }+\Omega _{r}\left( g,\delta
\right) _{p\left( \cdot \right) ,\omega }
\end{equation*}%
\begin{equation}
\leq \left( 1+24\mathbb{S}_{3}\mathbb{S}_{5}\mathbb{S}_{0}\right)
^{r}\left\Vert f\text{-}g\right\Vert _{p\left( \cdot \right) ,\omega }\text{+%
}12\mathbb{S}_{3}2^{r}\mathbb{S}_{5}\mathbb{S}_{0}\delta ^{r}\left\Vert
g^{\left( r\right) }\right\Vert _{p\left( \cdot \right) ,\omega }  \label{mn}
\end{equation}%
and taking infimum on $g\in W_{p\left( \cdot \right) ,\omega }^{r}$ in (\ref%
{mn}) we get%
\begin{equation*}
\Omega _{r}\left( f,\delta \right) _{p\left( \cdot \right) ,\omega }\leq
\left( 1+24\mathbb{S}_{3}2^{r}\mathbb{S}_{5}\mathbb{S}_{0}\right)
^{r}K_{r}\left( f,\delta ;p\left( \cdot \right) ,\omega \right) .
\end{equation*}
\end{proof}

\begin{proof}[\textbf{Proof of Theorem \protect\ref{TBir}}]
Proof of (\ref{b3}): Let $m\in \mathbb{N}$. Then%
\begin{equation*}
\left\Vert f^{\left( \alpha \right) }-D_{n,r}\left( f^{\left( \alpha \right)
}\right) \right\Vert _{p(\cdot ),\omega }=\left\Vert f^{\left( \alpha
\right) }\text{-}\mathcal{A}_{1/n}^{m}\left( f^{\left( \alpha \right)
}\right) \text{+}\mathcal{A}_{1/n}^{m}\left( f^{\left( \alpha \right)
}\right) -\right.
\end{equation*}%
\begin{equation*}
\left. -D_{n,r}\mathcal{A}_{1/n}^{m}\left( f^{\left( \alpha \right) }\right) 
\text{+}D_{n,r}\mathcal{A}_{1/n}^{m}\left( f^{\left( \alpha \right) }\right) 
\text{-}D_{n,r}\left( f^{\left( \alpha \right) }\right) \right\Vert
_{p(\cdot ),\omega }
\end{equation*}%
\begin{equation*}
\leq \left\Vert f^{\left( \alpha \right) }-\mathcal{A}_{1/n}^{m}\left(
f^{\left( \alpha \right) }\right) \right\Vert _{p(\cdot ),\omega
}+\left\Vert \mathcal{A}_{1/n}^{m}\left( f^{\left( \alpha \right) }\right)
-D_{n,r}\mathcal{A}_{1/n}^{m}\left( f^{\left( \alpha \right) }\right)
\right\Vert _{p(\cdot ),\omega }+
\end{equation*}%
\begin{equation*}
+\left\Vert D_{n,r}(\mathcal{A}_{1/n}^{m}\left( f^{\left( \alpha \right)
}\right) -f^{\left( \alpha \right) }\right\Vert _{p(\cdot ),\omega }\leq 
\mathbb{S}_{6}\Omega _{m}\left( f^{\left( \alpha \right) },\frac{1}{n}%
\right) _{p(\cdot ),\omega }+
\end{equation*}%
\begin{equation*}
+24\mathbb{S}_{3}2^{m}\mathbb{S}_{5}\mathbb{S}_{0}\left( \frac{2}{n^{m}}%
\left\Vert (\mathcal{A}_{1/n}^{r}f^{\left( \alpha \right) })^{\left(
m\right) }\right\Vert _{p(\cdot ),\omega }+\left\Vert \mathcal{A}%
_{1/n}^{m}\left( f^{\left( \alpha \right) }\right) -f^{\left( \alpha \right)
}\right\Vert _{p(\cdot ),\omega }\right)
\end{equation*}%
\begin{equation*}
\leq (\mathbb{S}_{6}+48\mathbb{S}_{3}2^{m}\mathbb{S}_{5}\mathbb{S}_{0}%
\mathbb{S}_{7}+24\mathbb{S}_{3}2^{m}\mathbb{S}_{5}\mathbb{S}_{0}\mathbb{S}%
_{6})\Omega _{m}\left( f^{\left( \alpha \right) },\frac{1}{n}\right)
_{p(\cdot ),\omega }
\end{equation*}%
\begin{equation*}
\leq \mathbb{S}_{8}K_{m}\left( f^{\left( \alpha \right) },n^{-1},p(\cdot
),\omega \right)
\end{equation*}%
with $\mathbb{S}_{8}:=(\mathbb{S}_{6}+48\mathbb{S}_{3}2^{m}\mathbb{S}_{5}%
\mathbb{S}_{0}\mathbb{S}_{7}+24\mathbb{S}_{3}2^{m}\mathbb{S}_{5}\mathbb{S}%
_{0}\mathbb{S}_{6})\left( 1+24\mathbb{S}_{3}2^{r}\mathbb{S}_{5}\mathbb{S}%
_{0}\right) ^{m}.$

Proof of (\ref{Eqn12}): For $g\in W_{p\left( \cdot \right) ,\omega }^{m}$ we
have%
\begin{equation*}
E_{n}\left( f^{(\alpha )}\right) _{p\left( \cdot \right) ,\omega }\leq \frac{%
48\mathbb{S}_{3}2^{r}\mathbb{S}_{5}\mathbb{S}_{0}}{n^{r}}E_{n}\left( f^{%
\text{ }(\alpha +r)}\right) _{p\left( \cdot \right) ,\omega }
\end{equation*}%
\begin{equation*}
\leq \frac{48\mathbb{S}_{3}2^{r}\mathbb{S}_{5}\mathbb{S}_{0}}{n^{r}}\left(
E_{n}\left( f^{\text{ }(\alpha +r)}-g\right) _{p\left( \cdot \right) ,\omega
}+E_{n}\left( g\right) _{p\left( \cdot \right) ,\omega }\right)
\end{equation*}%
\begin{equation*}
\leq \frac{48\mathbb{S}_{3}2^{r}\mathbb{S}_{5}\mathbb{S}_{0}}{n^{r}}\left[
\left\Vert f^{\text{ }(\alpha +r)}-g\right\Vert _{p\left( \cdot \right)
,\omega }+n^{-m}\left\Vert g^{\left( m\right) }\right\Vert _{p\left( \cdot
\right) ,\omega }\right]
\end{equation*}%
Taking infimum on $g\in W_{p\left( \cdot \right) ,\omega }^{m}$ in the last
inequality we find%
\begin{equation*}
E_{n}\left( f^{\text{ }(\alpha )}\right) _{p\left( \cdot \right) ,\omega
}n^{r}\leq 48\mathbb{S}_{3}2^{r}\mathbb{S}_{5}\mathbb{S}_{0}K_{m}\left( f^{%
\text{ }(\alpha +r)},n^{-1},p\left( \cdot \right) ,\omega \right)
\end{equation*}%
and inequality (\ref{Eqn12}) is proved.

Proof of (\ref{IJ1}): For $\rho \leq n,$

\begin{equation*}
\dfrac{K_{r}\left( f^{(\beta )},\rho ^{-1},p\left( \cdot \right) ,\omega
\right) }{\rho ^{\beta -\alpha }}\leq \dfrac{1}{\rho ^{\beta -\alpha }}\left[
1+\dfrac{n}{\rho }\right] ^{r}K_{r}\left( f^{(\beta )},n^{-1},p\left( \cdot
\right) ,\omega \right)
\end{equation*}%
and hence%
\begin{equation*}
\prod\limits_{\rho =1}^{n}\dfrac{K_{r}(f^{(\beta )},\rho ^{-1},p\left( \cdot
\right) ,\omega )}{\rho ^{\beta -\alpha }}\leq \prod\limits_{\rho =1}^{n}%
\frac{1}{\rho ^{\beta -\alpha }}\left[ 1+\frac{n}{\rho }\right]
^{r}K_{r}\left( f^{(\beta )},n^{-1},p\left( \cdot \right) ,\omega \right)
\end{equation*}%
\begin{equation*}
=\left( \frac{K_{r}\left( f^{\text{ }(\beta )},n^{-1},p\left( \cdot \right)
,\omega \right) }{n^{\beta -\alpha }}\right) ^{n}\prod_{\rho
=1}^{n}\left( 1+\frac{n}{\rho }\right) ^{r+\beta -\alpha }.
\end{equation*}%
Using Stirling's formula we have%
\begin{equation*}
\prod_{\rho =1}^{n}\left( 1+\frac{n}{\rho }\right) ^{r+\beta -\alpha
}\leq (2e)^{r+\beta -\alpha }
\end{equation*}%
and consequently%
\begin{equation*}
\left( \prod_{\rho =1}^{n}\frac{1}{\rho ^{\beta -\alpha }}K_{r}\left(
f^{\text{ }(\beta )},\rho ^{-1},p\left( \cdot \right) ,\omega \right)
\right) ^{1/n}\leq
\end{equation*}%
\begin{equation*}
\leq \left( 2e\right) ^{r+\beta -\alpha }\frac{1}{n^{\beta -\alpha }}%
K_{r}\left( f^{\text{ }(\beta )},n^{-1},p\left( \cdot \right) ,\omega
\right) .
\end{equation*}%
From (\ref{Eqn12}) and the property $E_{n}\left( f\right) _{p\left( \cdot
\right) }\downarrow $ (as $n\uparrow $), we find%
\begin{equation*}
\left( \prod_{\rho =1}^{n}E_{\rho }\left( f^{\text{ }(\alpha
)}\right) _{p\left( \cdot \right) ,\omega }\right) ^{1/n}\leq \left(
\prod\limits_{\rho =1}^{n}\frac{48\mathbb{S}_{3}2^{r}\mathbb{S}_{5}\mathbb{S}%
_{0}}{\rho ^{\beta -\alpha }}K_{r}\left( f^{\text{ }(\beta )},\rho
^{-1},p\left( \cdot \right) ,\omega \right) \right) ^{1/n}
\end{equation*}%
\begin{equation*}
\leq 48\mathbb{S}_{3}2^{r}\mathbb{S}_{5}\mathbb{S}_{0}\left( 2e\right)
^{r+\beta -\alpha }\frac{1}{n^{\beta -\alpha }}K_{r}\left( f^{\text{ }(\beta
)},n^{-1},p\left( \cdot \right) ,\omega \right) .
\end{equation*}%
and the result (\ref{IJ1}) follows.
\end{proof}

\begin{proof}[\textbf{Proof of Theorem \protect\ref{TTers}}]
Standard proof suits (see for example \cite{akg11gmj}) where used%
\begin{equation*}
\Omega _{r}\left( g,h\right) _{p\left( \cdot \right) ,\omega }\leq (12)^{r}%
\mathbb{S}_{3}^{r}\mathbb{S}_{5}^{r}\mathbb{S}_{0}^{r}h^{r}\left\Vert
g^{\left( r\right) }\right\Vert _{p\left( \cdot \right) ,\omega }\text{%
,\quad }\forall g\in W_{p\left( \cdot \right) ,\omega }^{r},
\end{equation*}%
and Bernstein type inequality%
\begin{equation*}
\left\Vert T^{(r)}\right\Vert _{p\left( \cdot \right) ,\omega }\leq 24%
\mathbb{S}_{3}\mathbb{S}_{5}\mathbb{S}_{0}2^{-r}n^{r}\left\Vert \Delta _{\pi
/n}^{r}T\right\Vert _{p\left( \cdot \right) ,\omega }\text{,\quad }\forall
T\in \mathcal{T}_{n}.
\end{equation*}
\end{proof}

\begin{proof}[\textbf{Proof of Theorem \protect\ref{Ters}}]
\begin{equation*}
\left\Vert \mathcal{U}_{\left( I-T_{t}\right) ^{r}f,F}\right\Vert _{C\left( 
\text{\textsc{T}}\right) }=\left\Vert \left( I-T_{t}\right) ^{r}\mathcal{U}%
_{f,F}\right\Vert _{C\left( \text{\textsc{T}}\right) }\leq \left\Vert \left(
I-\tilde{T}_{t}\right) ^{r}\mathcal{U}_{f,F}\right\Vert _{C\left( \text{%
\textsc{T}}\right) }
\end{equation*}%
\begin{equation*}
\leq 10\pi \left( 1+2^{2r-1}\right) 2^{2r+3k}\mathfrak{S}_{r+k}t^{r}\int%
\limits_{t}^{1}\frac{\omega _{r+k}\left( \mathcal{U}_{f,F},t\right)
_{C\left( \text{\textsc{T}}\right) }}{u^{r+1}}du
\end{equation*}%
\begin{equation*}
=10\pi \left( 1+2^{2r-1}\right) 2^{2r+3k}\mathfrak{S}_{r+k}2^{r}\mathfrak{S}%
_{r}t^{r}\int\limits_{t}^{1}\frac{\left\Vert \left( I-T_{t}\right) ^{r+k}%
\mathcal{U}_{f,F}\right\Vert _{C\left( \text{\textsc{T}}\right) }}{u^{r+1}}du
\end{equation*}%
\begin{equation*}
=\mathbb{S}_{9}10\pi \left( 1+2^{2r-1}\right) 2^{2r+3k}\mathfrak{S}%
_{r+k}2^{r}\mathfrak{S}_{r}t^{r}\int\limits_{t}^{1}\frac{\left\Vert \mathcal{%
U}_{\left( I-T_{t}\right) ^{r+k}f,F}\right\Vert _{C\left( \text{\textsc{T}}%
\right) }}{u^{r+1}}du.
\end{equation*}%
with $\mathbb{S}_{9}:=10\pi \left( 1+2^{2r-1}\right) 2^{2r+3k}\mathfrak{S}%
_{r+k}2^{r}\mathfrak{S}_{r}$.

Using TR we get%
\begin{equation*}
\frac{\left\Vert \left( I-T_{t}\right) ^{r}f\right\Vert _{p\left( \cdot
\right) ,\omega }}{t^{r}\int\limits_{t}^{1}\frac{\left\Vert \left(
I-T_{t}\right) ^{r+k}f\right\Vert _{p\left( \cdot \right) ,\omega }}{u^{r+1}}%
du}\leq 24\mathbb{S}_{3}\mathbb{S}_{5}\mathbb{S}_{0}\mathbb{S}_{9}.
\end{equation*}%
Then%
\begin{equation*}
\frac{K_{r}\left( f,t,p\left( \cdot \right) ,\omega \right) }{%
t^{r}\int\limits_{t}^{1}\frac{K_{r+k}\left( f,u,p\left( \cdot \right)
,\omega \right) }{u^{r+1}}du}\leq \text{max}\left\{ \mathbb{S}_{6},\mathbb{S}%
_{7}\right\} 24\mathbb{S}_{3}\mathbb{S}_{5}\mathbb{S}_{0}\mathbb{S}%
_{9}\left( 1\text{+}24\mathbb{S}_{3}2^{r}\mathbb{S}_{5}\mathbb{S}_{0}\right)
^{r+k}
\end{equation*}

as desired.
\end{proof}

\begin{proof}[\textbf{Proof of Theorem \protect\ref{simlt}}]
For $n\in \mathbb{N}_{0}$, we define dela Valee Poussin mean%
\begin{equation*}
W_{n}(f):=W_{n}(\cdot ,f):=\frac{1}{n+1}\sum\limits_{\nu =n}^{2n}S_{\nu
}(\cdot ,f)\in \mathcal{T}_{2n}\text{.}
\end{equation*}

It is well known that%
\begin{equation*}
W_{n}(f^{(\alpha )})=\left( W_{n}(f)\right) ^{(\alpha )}.
\end{equation*}

Proof of (\ref{fc}): Suppose that $q\in \mathcal{T}_{n}$ and $E_{n}\left(
f^{\left( k\right) }\right) _{p\left( \cdot \right) ,\omega }=\left\Vert
f^{\left( k\right) }-q\right\Vert _{p\left( \cdot \right) ,\omega }.$ Since 
\begin{equation*}
\mathcal{U}_{W_{n}\left( g\right) ,F}=W_{n}\left( \mathcal{U}_{g,F}\right)
\end{equation*}%
using TR we can find that%
\begin{equation*}
\left\Vert W_{n}\left( g\right) \right\Vert _{p\left( \cdot \right) ,\omega
}\leq 72\mathbb{S}_{3}\mathbb{S}_{5}\mathbb{S}_{0}\left\Vert g\right\Vert
_{p\left( \cdot \right) ,\omega }\quad \forall g\in L_{2\pi ,\omega
}^{p\left( \cdot \right) }.
\end{equation*}%
On the other hand we know that%
\begin{equation*}
W_{n}\left( t_{n}^{\ast }\right) =t_{n}^{\ast }\left( \cdot \right) .
\end{equation*}%
As a consequence%
\begin{equation*}
\left\Vert f^{\left( k\right) }\text{-}(t_{n}^{\ast })^{\left( k\right)
}\right\Vert _{p\left( \cdot \right) ,\omega }\leq \left\Vert f^{\left(
k\right) }\text{-}\left( W_{n}(f)\right) ^{\left( k\right) }\right\Vert
_{p\left( \cdot \right) ,\omega }\text{+}\left\Vert \left( W_{n}(f)\right)
^{\left( k\right) }\text{-}(t_{n}^{\ast })^{\left( k\right) }\right\Vert
_{p\left( \cdot \right) ,\omega }
\end{equation*}%
\begin{eqnarray*}
&\leq &\left\Vert f^{\left( k\right) }\text{-}q\right\Vert _{p\left( \cdot
\right) ,\omega }+\left\Vert q\text{-}W_{n}\left( f^{\left( k\right)
}\right) \right\Vert _{p\left( \cdot \right) ,\omega }+\left\Vert \left(
W_{n}(f)\text{-}t_{n}^{\ast }\right) ^{\left( k\right) }\right\Vert
_{p\left( \cdot \right) ,\omega } \\
&\leq &E_{n}\left( f^{\left( k\right) }\right) _{p\left( \cdot \right)
,\omega }\text{+}\left\Vert W_{n}\left( q\text{-}f^{\left( k\right) }\right)
\right\Vert _{p\left( \cdot \right) ,\omega }\text{+}24\mathbb{S}_{3}\mathbb{%
S}_{5}\mathbb{S}_{0}\mathfrak{S}_{k}2^{-k}(2n)^{k}\left\Vert W_{n}(f)\text{-}%
t_{n}^{\ast }\right\Vert _{p\left( \cdot \right) ,\omega } \\
&\leq &(1+72\mathbb{S}_{3}\mathbb{S}_{5}\mathbb{S}_{0})E_{n}\left( f^{\left(
k\right) }\right) _{p\left( \cdot \right) ,\omega }+24\mathbb{S}_{3}\mathbb{S%
}_{5}\mathbb{S}_{0}\mathfrak{S}_{k}n^{k}\left\Vert W_{n}(f\text{-}%
t_{n}^{\ast })\right\Vert _{p\left( \cdot \right) ,\omega } \\
&\leq &\frac{(1+72\mathbb{S}_{3}\mathbb{S}_{5}\mathbb{S}_{0})48\mathbb{S}_{3}%
\mathbb{S}_{5}\mathbb{S}_{0}2^{r}}{n^{r-k}}E_{n}\left( f^{\left( r\right)
}\right) _{p\left( \cdot \right) ,\omega }+ \\
&&+72\mathbb{S}_{3}\mathbb{S}_{5}\mathbb{S}_{0}48\mathbb{S}_{3}\mathbb{S}_{5}%
\mathbb{S}_{0}2^{r}n^{k}\frac{1}{n^{r}}E_{n}\left( f^{\left( r\right)
}\right) _{p\left( \cdot \right) ,\omega } \\
&=&\frac{\mathbb{S}_{10}}{n^{r-k}}E_{n}\left( f^{\left( r\right) }\right)
_{p\left( \cdot \right) ,\omega }
\end{eqnarray*}%
with $\mathbb{S}_{10}:=(1+72\mathbb{S}_{3}\mathbb{S}_{5}\mathbb{S}_{0})48%
\mathbb{S}_{3}\mathbb{S}_{5}\mathbb{S}_{0}2^{r}+72\mathbb{S}_{3}\mathbb{S}%
_{5}\mathbb{S}_{0}48\mathbb{S}_{3}\mathbb{S}_{5}\mathbb{S}_{0}2^{r}.$

Proof of (\ref{fc}): Using%
\begin{equation*}
\left\Vert f-W_{n}\left( f\right) \right\Vert _{p\left( \cdot \right)
,\omega }\leq \left\Vert f-t_{n}^{\ast }+t_{n}^{\ast }-W_{n}\left( f\right)
\right\Vert _{p\left( \cdot \right) ,\omega }
\end{equation*}%
\begin{equation*}
=\left\Vert f-t_{n}^{\ast }+W_{n}\left( t_{n}^{\ast }\right) -W_{n}\left(
f\right) \right\Vert _{p\left( \cdot \right) ,\omega }\leq E_{n}\left(
f\right) _{p\left( \cdot \right) ,\omega }\text{+}72\mathbb{S}_{3}\mathbb{S}%
_{5}\mathbb{S}_{0}\left\Vert f\text{-}t_{n}^{\ast }\right\Vert _{p\left(
\cdot \right) ,\omega }
\end{equation*}%
\begin{equation*}
=\left( 1+72\mathbb{S}_{3}\mathbb{S}_{5}\mathbb{S}_{0}\right) E_{n}\left(
f\right) _{p\left( \cdot \right) ,\omega }
\end{equation*}%
we obtain%
\begin{equation*}
\left\Vert f-W_{n}\left( f\right) \right\Vert _{p\left( \cdot \right)
,\omega }\leq \left( 1+72\mathbb{S}_{3}\mathbb{S}_{5}\mathbb{S}_{0}\right)
E_{n}\left( f\right) _{p\left( \cdot \right) ,\omega }
\end{equation*}%
\begin{equation*}
\leq \left( 1+72\mathbb{S}_{3}\mathbb{S}_{5}\mathbb{S}_{0}\right) 48\mathbb{S%
}_{3}\mathbb{S}_{5}\mathbb{S}_{0}2^{r}\text{max}\left\{ \mathbb{S}_{6},%
\mathbb{S}_{7}\right\} \frac{1}{n^{r}}\Omega _{s}\left( f^{\left( r\right) },%
\frac{1}{n}\right) _{p\left( \cdot \right) ,\omega }
\end{equation*}%
\begin{equation*}
=\mathbb{S}_{11}\frac{1}{n^{r}}\Omega _{s}\left( f^{\left( r\right) },\frac{1%
}{n}\right) _{p\left( \cdot \right) ,\omega }
\end{equation*}%
with $\mathbb{S}_{11}:=\left( 1+72\mathbb{S}_{3}\mathbb{S}_{5}\mathbb{S}%
_{0}\right) 48\mathbb{S}_{3}\mathbb{S}_{5}\mathbb{S}_{0}2^{r}$max$\left\{ 
\mathbb{S}_{6},\mathbb{S}_{7}\right\} $.

From%
\begin{equation*}
\left\Vert f-t_{n}^{\ast }\right\Vert _{p\left( \cdot \right) ,\omega }\leq 
\frac{48\mathbb{S}_{3}\mathbb{S}_{5}\mathbb{S}_{0}2^{r}\text{max}\left\{ 
\mathbb{S}_{6},\mathbb{S}_{7}\right\} }{n^{r}}\Omega _{s}\left( f^{\left(
r\right) },\frac{1}{n}\right) _{p\left( \cdot \right) ,\omega }
\end{equation*}%
we obtain%
\begin{equation*}
\left\Vert W_{n}\left( f\right) -t_{n}^{\ast }\right\Vert _{p\left( \cdot
\right) ,\omega }\leq \frac{\mathbb{S}_{12}}{n^{r}}\Omega _{s}\left(
f^{\left( r\right) },\frac{1}{n}\right) _{p\left( \cdot \right) ,\omega }
\end{equation*}%
with $\mathbb{S}_{12}:=\mathbb{S}_{11}+48\mathbb{S}_{3}\mathbb{S}_{5}\mathbb{%
S}_{0}2^{r}$max$\left\{ \mathbb{S}_{6},\mathbb{S}_{7}\right\} .$ Hence%
\begin{equation*}
\left\Vert f^{\left( k\right) }\text{-}\left( W_{n}\left( f\right) \right)
^{\left( k\right) }\right\Vert _{p\left( \cdot \right) ,\omega }\leq
\left\Vert f^{\left( k\right) }\text{-}\left( t_{n}^{\ast }\right) ^{\left(
k\right) }\right\Vert _{p\left( \cdot \right) ,\omega }+\left\Vert \left(
W_{n}\left( f\right) \right) ^{\left( k\right) }\text{-}\left( t_{n}^{\ast
}\right) ^{\left( k\right) }\right\Vert _{p\left( \cdot \right) ,\omega }
\end{equation*}%
\begin{equation*}
\leq \mathbb{S}_{10}n^{k-r}E_{n}\left( f^{\left( r\right) }\right) _{p\left(
\cdot \right) ,\omega }+24\mathbb{S}_{12}\mathbb{S}_{3}\mathbb{S}_{5}\mathbb{%
S}_{0}\mathfrak{S}_{k}n^{k}\frac{1}{n^{r}}\Omega _{s}\left( f^{\left(
r\right) },\frac{1}{n}\right) _{p\left( \cdot \right) ,\omega }
\end{equation*}%
\begin{equation*}
\leq \mathbb{S}_{13}n^{k-r}\Omega _{s}\left( f^{\left( r\right) },1/n\right)
_{p\left( \cdot \right) ,\omega }
\end{equation*}%
with $\mathbb{S}_{13}:=\mathbb{S}_{10}+24\mathbb{S}_{12}\mathbb{S}_{3}%
\mathbb{S}_{5}\mathbb{S}_{0}\mathfrak{S}_{k}$.
\end{proof}

\begin{proof}[\textbf{Proof of Lemma \protect\ref{bukun}}]
Let $0<h\leq \delta <\infty $, $\mathcal{P}^{\log }$, $\omega \in A_{p\left(
\cdot \right) }$ and $f\in L_{2\pi ,\omega }^{p\left( \cdot \right) }$.
Then, 
\begin{equation*}
\left\Vert \mathcal{U}_{\left( I-T_{h}\right) f,F}\right\Vert _{C\left( 
\text{\textsc{T}}\right) }=\left\Vert \left( I-T_{h}\right) \mathcal{U}%
_{f,F}\right\Vert _{C\left( \text{\textsc{T}}\right) }\leq 72\left\Vert
\left( I-T_{\delta }\right) \mathcal{U}_{f,F}\right\Vert _{C\left( \text{%
\textsc{T}}\right) }
\end{equation*}%
\begin{equation*}
\leq 72\left\Vert \mathcal{U}_{\left( I-T_{\delta }\right) f,F}\right\Vert
_{C\left( \text{\textsc{T}}\right) }.
\end{equation*}%
From TR we get%
\begin{equation*}
\left\Vert \left( I-T_{h}\right) f\right\Vert _{p\left( \cdot \right)
,\omega }\leq 1728\mathbb{S}_{3}\mathbb{S}_{5}\mathbb{S}_{0}\left\Vert
\left( I-T_{\delta }\right) f\right\Vert _{p\left( \cdot \right) ,\omega }.
\end{equation*}
\end{proof}

\begin{proof}[\textbf{Proof of Lemma \protect\ref{bukunA}}]
If $\mathcal{P}^{\log }$, $\omega \in A_{p\left( \cdot \right) }$ and $f\in
L_{2\pi ,\omega }^{p\left( \cdot \right) }$, then, using generalized
Minkowski's integral inequality and Lemma \ref{bukun} we obtain%
\begin{equation*}
\left\Vert \left( I-\mathfrak{R}_{\delta }\right) f\right\Vert _{p\left(
\cdot \right) }=\left\Vert \frac{2}{\delta }\int\nolimits_{\delta
/2}^{\delta }\left( \frac{1}{h}\int\nolimits_{0}^{h}\left( f\left(
x+t\right) -f\left( x\right) \right) dt\right) dh\right\Vert _{p\left( \cdot
\right) }
\end{equation*}%
\begin{equation*}
=\left\Vert \frac{2}{\delta }\int\nolimits_{\delta /2}^{\delta }\left(
T_{h}f\left( x\right) -f\left( x\right) \right) dh\right\Vert _{p\left(
\cdot \right) }\leq \frac{2}{\delta }\int\nolimits_{\delta /2}^{\delta
}\left\Vert T_{\delta }f-f\right\Vert _{p\left( \cdot \right) }dh
\end{equation*}%
\begin{equation*}
\leq 1728\mathbb{S}_{3}\mathbb{S}_{5}\mathbb{S}_{0}\left\Vert T_{\delta
}f-f\right\Vert _{p\left( \cdot \right) }\frac{2}{\delta }%
\int\nolimits_{\delta /2}^{\delta }dh=1728\mathbb{S}_{3}\mathbb{S}_{5}%
\mathbb{S}_{0}\left\Vert \left( I-T_{\delta }\right) f\right\Vert _{p\left(
\cdot \right) }.
\end{equation*}
\end{proof}

\begin{proof}[\textbf{Proof of Lemma \protect\ref{lm05}}]
Using%
\begin{equation*}
\left\Vert \mathcal{U}_{\delta \left( \mathfrak{R}_{\delta }f\right)
^{\prime },F}\right\Vert _{C\left( \text{\textsc{T}}\right) }=\left\Vert
\delta \left( \mathcal{U}_{\left( \mathfrak{R}_{\delta }f\right) ,F}\right)
^{\prime }\right\Vert _{C\left( \text{\textsc{T}}\right) }=\delta \left\Vert
\left( \mathfrak{R}_{\delta }(\mathcal{U}_{f,F})\right) ^{\prime
}\right\Vert _{C\left( \text{\textsc{T}}\right) }
\end{equation*}%
\begin{equation*}
\leq \cdots \leq 2\left( 37+146\ln 2^{36}\right) \left\Vert \left(
I-T_{\delta }\right) (\mathcal{U}_{f,F})\right\Vert _{C\left( \text{\textsc{T%
}}\right) }
\end{equation*}%
\begin{equation*}
=2\left( 37+146\ln 2^{36}\right) \left\Vert (\mathcal{U}_{\left( I-T_{\delta
}\right) f,F})\right\Vert _{C\left( \text{\textsc{T}}\right) }
\end{equation*}%
we conclude from TR that%
\begin{equation*}
\delta \left\Vert (\mathfrak{R}_{\delta }f)^{\prime }\right\Vert _{p\left(
\cdot \right) }\leq \mathbb{S}_{3}\mathbb{S}_{5}\mathbb{S}_{0}48\left(
37+146\ln 2^{36}\right) \left\Vert \left( I-T_{\delta }\right) f\right\Vert
_{p\left( \cdot \right) }.
\end{equation*}
\end{proof}


\begin{thebibliography}{99}
\bibitem{ascs} F. Abdullaev, A. Shidlich and S. Chaichenko, \textit{Direct
and inverse approximation theorems of functions in the Orlicz type spaces},
Math. Slovaca, 69 (2019), 1367-1380.

\bibitem{aosss} F. Abdullaev, N. \"{O}zkaratepe, V. Savchuk and A. Shidlich, 
\textit{Exact constants in direct and inverse approximation theorems for
functions of several variables in the spaces }$S_{p}$, Filomat, 33 (2019),
1471-1484.

\bibitem{ahh} T. Adamowicz, P. Harjulehto and P. H\"{a}st\"{o}, \textit{%
Maximal Operator in Variable Exponent Lebesgue Spaces on Unbounded
Quasimetric Measure Spaces}, Mathematica Scandinavica, 116 (2015), no: 1,
5-22.

\bibitem{Akgeja} R. Akg\"{u}n, \textit{Polynomial approximation in weighted
Lebesgue spaces}, East J. Approx., 17 (2011), no. 3, 253-266.

\bibitem{ra} R. Akg\"{u}n, \textit{Approximating polynomials for functions
of weighted Smirnov-Orlicz spaces}, J. Funct. Spaces Appl., Volume 2012,
Article ID 982360, 41 pages.

\bibitem{ra19} R. Akg\"{u}n, \textit{Direct theorems of trigonometric
approximation for variable exponent Lebesgue spaces}, Revista de la Uni\'{o}%
n Matem\'{a}tica Argentina, 60 (2019), no. 1, 121-135.

\bibitem{ra-dmi} ---, \textit{Simultaneous and converse approximation
theorems in weighted Orlicz spaces}, Bull. Belg. Math. Soc. Simon Stevin, 
\textbf{17} (2010), no:1, pp. 13-28.

\bibitem{ra11u} R. Akg\"{u}n, \textit{Trigonometric approximation of
functions in generalized Lebesgue spaces with variable exponent}, Ukrainian
Math. J., 63 (2011), no. 1, 1-26.

\bibitem{akg11gmj} R. Akg\"{u}n, \textit{Polynomial approximation of
functions in weighted Lebesgue and Smirnov spaces with non-standart growth},
Georgian Math. J., 18 (2011), no. 2, 203-235.

\bibitem{AK1} R. Akg\"{u}n, \textit{Sharp Jackson and converse theorems of
trigonometric approximation in weighted Lebesgue spaces}, Proc. A. Razmadze
Math. Inst., 152 (2010), pp. 1-18.

\bibitem{AG} R. Akg\"{u}n; A. Ghorbanalizadeh, \textit{Approximation by
integral functions of finite degree in variable exponent Lebesgue spaces on
the real axis}. Turk. J. Math. 42 (2018), no. 4, 1887-1903.

\bibitem{ahak} A.H. Av\c{s}ar and H. Ko\c{c}, \textit{Jackson and Stechkin
type inequalities of trigonometric approximation in }$A_{p,q(.)}^{w,\theta }$%
, Turk J Math 42 (2018), no: , 2979-2993.

\bibitem{ay12} I. Aydin, Weighted Variable Sobolev Spaces and Capacity,
Journal of Function Spaces and Applications, Volume 2012, Article ID 132690,
17 pages.

\bibitem{au20} I. Aydin and C. Unal, \textit{The Kolmogorov--Riesz theorem
and some compactness criterions of bounded subsets in weighted variable
exponent amalgam and Sobolev spaces}, Collectanea Mathematica (2020) 71,
no:3, 349-367.

\bibitem{B12} R. A. Bandaliev, \textit{Application of multimensional Hardy
operator and its connection with a certain nonlinear differential equation
in weighted variable Lebesgue spaces}, Ann. Funct. Anal., 4 (2013), no. 2,
118-130.

\bibitem{bg03} E. Berkson and T. A. Gillespie, \textit{On restrictions of
multipliers in weighted setting}, Indian U. Math. J., 52 (2003), no:4,
927-961.

\bibitem{ch12} V. Chaichenko, \textit{Best approximation of periodic
functions in generalized Lebesgue spaces}, Ukrainian Math. J., 64, (2013),
no. 9, 1421-1439.

\bibitem{cuf} D. Cruz-Uribe SFO, A. Fiorenza, \textit{Approximate identities
in variable }$Lp$\textit{\ spaces}. Mathematische Nachrichten 2007; 280 (3):
256-270.

\bibitem{UF13} D. Cruz-Uribe SFO, A. Fiorenza, Variable Lebesgue Spaces,
Foundations and Harmonic Analysis, Birkhauser, Applied and Numerical
Harmonic Analysis, 2013.

\bibitem{cufb} D. Cruz-Uribe SFO, A. Fiorenza, C. J. Neugebauer, \textit{%
Weighted norm inequalities for the maximal operator on variable Lebesgue
spaces}, Journal of mathematical analysis and applications, (2012), 394,
744-760.

\bibitem{cudh} D. Cruz-Uribe SFO, L. Diening, P. H\"{a}st\"{o}, \textit{The
maximal operator on weighted variable Lebesgue spaces}, Fract. Calc. Appl.
Anal. (2011), 14, 361-374.

\bibitem{devore} R. A. Devore, G. G. Lorentz, Constructive Approximation,
Springer-Verlag, (1993).

\bibitem{DHHR11} L. Diening, P. Harjulehto, P. H\"{a}st\"{o} and M. Ru\v{z}i%
\v{c}ka, Lebesgue and Sobolev spaces with variable exponents, Lecture Notes
in Math., vol. 2017, Springer, Berlin, Heidelberg, 2011.

\bibitem{d-h} L. Diening and P. H\"{a}st\"{o}, \textit{Muckenhoupt weights
in variable exponent spaces}, Albert Ludwings Universit\"{a}t Freiburg,
Mathematische Fakult\"{a}t, https://www.problemsolving.fi/pp/p75\_submit.pdf

\bibitem{z76} Z. Ditzian, \textit{Inverse theorems for functions in }$L^{p}$%
\textit{\ and other spaces}, Proc. Amer. Math. Soc. 54 (1976), 80-82.

\bibitem{day} A. Dogu, A.H. Avsar and Y.E. Yildirir, \textit{Some
inequalities about convolution and trigonometric approximation in weighted
Orlicz spaces}, Proceedings of the Institute of Mathematics and Mechanics,
National Academy of Sciences of Azerbaijan, Volume 44, Number 1, 2018, Pages
107-115.

\bibitem{douux} J. Duoandikoetxea, Fourier Analysis, Graduate Studies in
Mathematics Volume 29, AMS, Providence, Rhode Island, 2000.

\bibitem{DzSh} V. K. Dzyadyk and I. A. Shevchuk, Theory of Uniform
Approximation of Functions by Polynomials, De Gruyter, 2008.

\bibitem{eah86} E. A. Gadjieva, \textit{Investigation of the properties of
functions with quasimonotone Fourier coefficients in generalized
Nikolskii-Besov spaces}, author's summary of dissertation, Tbilisi, 1986,
(In Russian).

\bibitem{garfra85} J. Garcia-Cuerva, J. L. R. De Francia, Weighted Norm
Inequalities and Related Topics, Volume 116, Pages ii-viii, 1-604 (1985).

\bibitem{git20} D. V. Gorbachev V. I. Ivanov, \textit{Fractional smoothness
in }$L_{p}$\textit{\ with Dunkl weight and its applications}, Mathematical
Notes, 106 (2019), no:4, 537-561.

\bibitem{graMod} L. Grafakos, Modern Fourier Analysis, 2nd Ed., Springer,
2009.

\bibitem{HH} P. Harjulehto and P. H\"{a}st\"{o}, Orlicz Spaces and
Generalized Orlicz Spaces, Lecture Notes in Mathematics 2236, 2019.

\bibitem{dmi02} D. M. Israfilov, \textit{Approximation by p-Faber
polynomials in the weighted Smirnov class }$E^{p}$\textit{(}$G$\textit{,}$%
\omega $\textit{) and the Bieberbach polynomials}. Constr. Approx. 17
(2001), no. 3, 335-351.

\bibitem{isrgu10} D. M. Israfilov and A. Guven, \textit{Trigonometric
approximation in generalized Lebesgue spaces}$\ L^{p(x)}$, J. Math. Inequal.
4 (2010), no. 2, 285-299.

\bibitem{Israfil} D. M. Israfilov and A. Testici, \textit{Approximation
problems in the Lebesgue spaces with variable exponent}, Journal of
Mathematical Analysis and Applications, \textit{459}:1, 2018, 112--123.

\bibitem{isrYir16} D. M. Israfilov and E. Yirtici, \textit{Convolutions and
best approximations in variable exponent Lebesgue spaces}, Math. Reports,
18(68) (2016), no. 4, 497-508.

\bibitem{Ja} S.Z. Jafarov, \textit{Linear Methods for Summing Fourier Series
and Approximation in Weighted Lebesgue Spaces with Variable Exponents}.
Ukrainian Mathematical Journal. 2015; 66(10): 1509-1518.

\bibitem{Ja1} S.Z. Jafarov, \textit{Approximation by trigonometric
polynomials in subspace of variable exponent grand Lebesgue spaces}. Global
Journal of Mathematics, 2016; 8(2): 836--843.

\bibitem{kc16} H. Koc, \textit{Simultaneous approximation by polynomials in
Orlicz spaces generated by quasiconvex Young functions}, Kuwait J. Sci. 43
(2016), no. 4, 18-31.

\bibitem{kol17} Y. Kolomoitsev, \textit{On moduli of smoothness and averaged
differences of fractional order}, Fractional Calculus and Applied Analysis,
20 (2017), no. 4, 988-1009.

\bibitem{kyil10} V. Kokilashvili, and Y. E. Yildirir, \textit{On the
approximation by trigonometric polynomials in weighted Lorentz spaces}. J.
Funct. Spaces Appl. 8 (2010), no. 1, 67-86.

\bibitem{vmk-sgs03} V. Kokilashvili and S. G. Samko, \textit{Singular
integrals weighted Lebesgue spaces with variable exponent}, Georgian M. J., 
\textbf{10}-1 (2003), 145-156.

\bibitem{Ky97} N. X. Ky, \textit{Moduli of mean smoothness and approximation
with }$\ A_{p}$\textit{-weights}, Annales Univ. Sci. Budapest, Sectio Math.,
40 (1997), 37-48.

\bibitem{lsz20} W. \L enski, and B. Szal, \textit{Trigonometric
approximation of functions from }$L_{2\pi }^{p(x)}$. Results Math. 75
(2020), no. 2, Paper No. 56, 14 pp.

\bibitem{muc} B. Muckenhoupt, \textit{Weighted norm inequalities for the
Hardy maximal function}, Trans Amer Math Soc \textbf{165} (1972), 207-226.

\bibitem{TiNa} G. I. Natanson and M. F. Timan, \textit{The geometric means
of the sequence of best approximations}, (Russian) Vestnik Leningrad. Univ.
Mat. Mekh. Astronom., 1979, vyp. 4, 50-52.

\bibitem{Ni} S. M. Nikolskii, \textit{Inequalities for entire functions of
finite degree and their application to the theory of differentiable
functions of several variables}, Amer. Math. Soc. Transl. Ser. 2, 80 (1969),
1-38, (Trudy Mat. Inst. Steklov 38 (1951), 211-278).

\bibitem{sn21} Y. Sawano T. Nogayama, \textit{Local Muckenhoupt class for
variable exponents}, 2021, 2021:70, 1-27.

\bibitem{sh08} I. I. Sharapudinov, \textit{Some problems in approximation
theory in the spaces }$L^{p(x)}(E)$, (Russian) Anal. Math. 33 (2007), no. 2,
135-153.

\bibitem{Sh12} I. I. Sharapudinov, Some questions in the theory of
approximation in Lebesgue spaces with variable exponent, Itogi Nauki. Yug
Rossii. Mat. Monografiya, vol. 5, Vladikavkaz 2012, 267 pp. Russian.

\bibitem{sh13i} I. I. Sharapudinov, \textit{Approximation of functions in }$%
L_{2\pi }^{p(\cdot )}$\textit{\ by trigonometric polynomials}, Izv. RAN.
Ser. Mat., 77 (2013), no. 2, 197-224.

\bibitem{sh15} I. I. Sharapudinov, \textit{On direct and inverse theorems of
approximation theory in variable Lebesgue and Sobolev spaces}, Azerbaijan J.
Math., 4 (2014), no. 1, 55-72.

\bibitem{shEm14} T. N. Shakh-Emirov, \textit{On Uniform Boundedness of Some
Families of Integral Convolution Operators in Weighted Variable Exponent
Lebesgue Spaces}, Izv. Saratov. Univ. Mat. Mekh. Inform., 14 (2014), no.
4(1), 422-427.

\bibitem{stei} E. M. Stein, Harmonic analysis: real-variable methods,
orthogonality, and oscillatory integrals, Princeton University Press, 1993.

\bibitem{vol17} S. S. Volosivets, \textit{Approximation of functions and
their conjugates in variable Lebesgue spaces}, Sbornik: Mathematics, 208
(2017), no. 1, 44-59.

\bibitem{yeydmi10} Y. E. Yildirir and D. M. Israfilov, \textit{Approximation
theorems in weighted Lorentz spaces}, Carpathian J. Math., 26 (2010), No:1,
108-119.
\end{thebibliography}
\end{document}